\documentclass[12pt]{article}
\usepackage{amsmath,amssymb,amsthm}
\usepackage{color}
\numberwithin{equation}{section}
\newtheorem{thm}{Theorem}[section]
\newtheorem{prop}[thm]{Proposition}
\newtheorem{cor}[thm]{Corollary}
\newtheorem{exam}[thm]{Example}
\newtheorem{rem}[thm]{Remark}
\newtheorem{lem}[thm]{Lemma}
\newtheorem{defn}[thm]{Definition}

\newcommand{\pf}{{\em Proof.}}

\def\address#1#2{\begingroup
\noindent\parbox[t]{16cm}{%
\small{\scshape\ignorespaces#1}\par\vskip1ex
\noindent\small{\itshape E-mail address}%
\/: #2\par\vskip4ex}\hfill%
\endgroup}%

\begin{document}
\title{Spread rate of branching Brownian motions}
\author{Yuichi Shiozawa
\footnote{Dedicated to Professor Masayoshi Takeda on the occasion of his 60th birthday}
\footnote{Supported in part by JSPS KAKENHI No.\ JP26400135, JP17K05299.}}
\maketitle 

\begin{abstract}

We find the exponential growth rate of the population outside a ball 
with time dependent radius for a branching Brownian motion in Euclidean space.  
We then see that the upper bound of the particle range is determined by 
the principal eigenvalue of the Schr\"odinger type operator associated with 
the branching rate measure and branching mechanism. 
We assume that the branching rate measure is small enough at infinity,  
and can be singular with respect to the Lebesgue measure. 
We finally apply our results to several concrete models.
\end{abstract}

\section{Introduction}

In this paper, we are concerned with the forefront of particles 
for a branching Brownian motion in Euclidean space. 
We first discuss the exponential growth rate of the population 
outside an open ball with radius increasing linearly in time. 
We then see that the forefront is determined by the principal eigenvalue 
of the Schr\"odinger type operator associated 
with the branching rate measure and branching mechanism. 
We finally apply our results to several concrete models.
\medskip

The Galton-Watson (branching) process is a stochastic model 
describing the evolution of a population; 
particles reproduce at each generation according to a common offspring distribution. 
One of the main interests in this model is the population growth rate. 
By martingale theory,  we know that on the non-extinction event, 
the total population is eventually comparable to its expectation 
(see, e.g., \cite[p.9, Theorems 1, 2]{AN04}). 
This result is valid also for the continuous time version in which 
the particle lifetimes are independently and exponentially distributed  
(see, e.g., \cite[p.111, Theorem 1 and p.112, Theorem 2]{AN04}). 

In order to make the Galton-Watson process more applicable and interesting, 
we consider a stochastic model in which 
particles reproduce according to a Galton-Watson process 
and move as Brownian motions on ${\mathbb R}^d$. 
This model is called the branching Brownian motion on ${\mathbb R}^d$ 
(see Subsection \ref{subsect-model} for definition). 
The distribution of each particle's lifetime depends on the size of a measure $\mu$ on ${\mathbb R}^d$ 
along the trajectory of that particle; we explain this in more detail in Subsection \ref{subsect-model}. 
The offspring distribution is given by a probability function $\{p_n(x)\}_{n\geq 1}$; 
each particle splits into $n$ particles with probability $p_n(x)$ at branching cite $x\in {\mathbb R}^d$.
We call $\mu$ and $\{p_n(x)\}_{n\geq 1}$, respectively, 
the branching rate (measure) and branching mechanism of the model. 
If $\mu$ is absolutely continuous with respect to the $d$-dimensional Lebesgue measure, 
then the density function is called the branching rate function. 
For instance, if the branching rate function is a constant function taking a value $c>0$, 
then the particle lifetimes are exponentially distributed with rate $c$. 
We refer to \cite{INW68-1,INW68-2, INW69} for 
the construction and properties of more general branching Markov processes.

The population at large times is well described by the bottom of the spectrum of 
the Schr\"odinger type operators. 
Let $Q(x):=\sum_{n=1}^{\infty}np_n(x)$ be 
the expected offspring number at branching cite $x\in {\mathbb R}^d$ and 
define the (formal) Schr\"odinger type operator
$${\cal H}^{(Q-1)\mu}:=-\frac{1}{2}\Delta-(Q-1)\mu,$$
where $\Delta$ is the Laplace operator on $L^2({\mathbb R}^d)$.
Denote by $\lambda$ and $h$, respectively,  
the bottom of the spectrum of ${\cal H}^{(Q-1)\mu}$ 
and corresponding eigenfunction (if it exists). 
Then $-\lambda$ becomes large in proportion to the branching intensity 
determined by the branching rate and mechanism.
Furthermore, $\lambda$ and $h$ reflect the long time asymptotic properties of 
branching Brownian motions:
under some conditions on the branching rate and mechanism, 
$\lambda$ characterizes the exponential population growth rate  
and $h$ determines the spatial distribution of the particles.  
These facts are known for more general branching Markov processes 
(see, e.g., \cite{CRY17+,CS07,EHK10,S08,W67}).
\medskip

Here we discuss the upper bound of the particle range for all sufficiently large time.  
For $d=1$, let $R_t$ be the maximal displacement, that is, 
the position of the right most particle at time $t$. 
If the particle lifetimes are independently and  exponentially distributed with rate $c>0$, 
and if each particle splits into two at branching, 
then $R_t/t\rightarrow \sqrt{2c}$ as $t\rightarrow\infty$
(see, e.g., \cite{B78} and \cite{R13} for more detailed properties on $R_t$.
See also \cite{M15} for an extension of \cite{B78} to the multidimensional case). 
This reveals how the intensity of branching appears in the particle range.  
A similar phenomenon occurs 
when the branching rate is a continuous function  
taking  a positive constant value outside some compact set (see \cite[Section 5]{LS89} for details). 
Erickson \cite{E84} proved that even if the branching rate function is small at infinity, 
$R_t/t\rightarrow\sqrt{-\lambda/2}$ as $t\rightarrow\infty$.
He also remarked in \cite[Section 4]{E84} that 
\begin{enumerate}
\item
this result also holds 
if for the Schr\"odinger type operator associated with the branching rate measure and branching mechanism,  
the principal eigenfunction  
decays exponentially at infinity (for example, this condition is satisfied if the measure is compactly supported);
\item 
a similar result is true for the multidimensional case 
if the branching rate function and branching mechanism  
are spherically symmetric.   
\end{enumerate}
In connection with (i), 
when $d=1$ and the branching rate is given by the Dirac measure, 
Bocharov and Harris \cite{BH14} determined 
the growth rate of the number of particles on some time dependent interval. 
As a corollary, they obtained the same result as in \cite{E84}. 
We know that if the branching rate function is given by $c|x|^p$ ($c>0$, $p\geq 0$), 
then $R_t$ grows faster than linear for $p\in (0,2]$,
and the branching Brownian motion explodes for $p>2$ 
(see, e.g., \cite{BBHH10,E84,HH09,IM74}). 

Our purpose in this paper is to show that 
if the measure $(Q-1)\mu$ satisfies a certain Kato class condition and $\lambda<0$,  
then the next assertions hold 
(see Theorem \ref{thm-main} and Corollary \ref{cor-main} for details).
\begin{itemize}
\item For $\delta\geq 0$, let $Z_t^{\delta t}$ be the number of particles at time $t$ 
outside an open ball with radius $\delta t$ centered at the origin.  
If $\delta>\sqrt{-\lambda/2}$, then 
$$\lim_{t\rightarrow\infty}Z_t^{\delta t}=0, \quad \text{a.s.}$$
On the other hand, if $0\leq \delta<\sqrt{-\lambda/2}$, 
then on an event of positive probability on some associated martingale, 
$$\lim_{t\rightarrow\infty}\frac{1}{t}\log Z_t^{\delta t}
=-\lambda-\sqrt{-2\lambda}\delta.$$
\item Let $L_t$ be the maximum of the Euclidean norms of particles alive at time $t$.  
Then on an event of positive probability on some associated martingale, 
\begin{equation}\label{eq-l-lln}
\lim_{t\rightarrow\infty}\frac{L_t}{t}=\sqrt{-\frac{\lambda}{2}}.
\end{equation}
\end{itemize}
The first and the latter are extensions of results,  respectively, 
of \cite{BH14} and of \cite{BH14,E84} 
to multidimensional branching Brownian motions 
with more general branching rate. 
We emphasize that 
the branching rate measure is allowed to be spherically asymmetric 
or to be singular
with respect to the Lebesgue measure. 
For $d=1$ , \eqref{eq-l-lln} is true 
if $L_t$ is replaced by $R_t$ (Theorem \ref{thm-direction} and Remark \ref{rem-mtg}). 
Moreover, we give explicitly conditions on the branching rate  
for the growth rate of $R_t$ being clarified.   
Note that if $\mu$ is compactly supported, $Q$ is bounded and if  $\lambda=0$, 
then no branching occurs for all sufficiently large time and hence 
we get a result of the law of the iterated logarithm type 
(Proposition \ref{prop-zero}).

Our approach here is similar to that of Bocharov and Harris \cite{BH14}. 
To obtain  the upper bound of $Z_t^{\delta t}$ in Lemma \ref{lem-limsup}, 
we determine the growth rate of the expectation of $Z_t^{\delta t}$ 
(Propositions \ref{prop-asymp} and \ref{prop-non-cpt}). 
This expectation coincides with that of the Feynman-Kac functional 
under the condition that the Brownian particle at time $t$ is outside 
an open ball with radius $\delta t$ centered at the origin (see \eqref{eq-br-fk} below).  
When $d=1$ and the branching rate is the Dirac measure, 
the Feynman-Kac functional is the exponential of the local time. 
Bocharov and Harris \cite{BH14} then calculated 
the expectation of $Z_t^{\delta t}$ explicitly 
by using the joint distribution of the Brownian motion and local time.  
However, this approach is not applicable to branching Brownian motions with more general branching rate.
To overcome this difficulty,  
we make use of the $L^p$-independence of the spectral bound of ${\cal H}^{(Q-1)\mu}$ 
(see \eqref{eq-comp-funct} below), 
which was established by Takeda \cite{T98,T03,T08}; in particular,  
a key fact is the asymptotics of the expectation of the Feynman-Kac functional  
(see (\ref{eq-feyn-lim}) below).  
The price is to impose a certain Kato class condition on the measure $(Q-1)\mu$. 
This condition means that the measure $(Q-1)\mu$ is small enough at infinity. 
To obtain the lower bound of $Z_t^{\delta t}$ in Lemma \ref{lem-liminf},  
we modify the argument of \cite[Lemma 4]{BH14} for the multidimensional case   
(compare the proof of Lemma \ref{lem-liminf} below with that of \cite[Lemma 4]{BH14}).
\medskip

Carmona and Hu \cite{CH14} studied the growth rate of the maximal displacement
of a branching random walk on ${\mathbb Z}$ 
such that each of particles moves as a  general irreducible (non-symmetric) random walk and  
branching occurs only on finite points. 
Bulinskaya \cite{B16+} further developed this study to the multidimensional case,  
and verified in which direction 
there exists a particle taking the maximal norm. 
In contrast with \cite{B16+} and \cite{CH14}, 
if the measure $(Q-1)\mu$ satisfies a certain Kato class condition and  $\lambda<0$,  
then the spread rate is uniform in direction for our model (Theorem \ref{thm-direction})
because particles move as Brownian motions on ${\mathbb R}^d$;
however, the branching rate measure can have the non-compact support, 
and its density function (if it exists) can be unbounded on some compact set 
(see Subsection \ref{subsect-exam} for examples).
\medskip

As we see by \eqref{eq-l-lln},  $L_t$ grows linearly 
even if the measure $(Q-1)\mu$ is small enough at infinity.  
If we drop this condition, then \eqref{eq-l-lln} is not true in general; 
as we mentioned before, 
if $d=1$ and the branching rate function is $c|x|^p$ with $c>0$ and $p\in (0,2]$, 
then $R_t$ grows faster than linear. 
On the other hand, if $d=1$ and 
the branching rate function takes a positive constant value outside some compact set, 
then $R_t/t$ converges to a positive constant. 
Combining this with \eqref{eq-l-lln}, we guess that even for $d\geq 2$, 
we will have a result similar to \eqref{eq-l-lln}  
if the branching rate function satisfies a Kato class condition 
(see Definition \ref{def-kato} (i) and Example \ref{ex-v} (i)) 
and takes a positive constant value outside some compact set. 
\medskip

We finally remark the relation between our results  and semilinear heat equations, 
which was already discussed in \cite{E84}. 
Let ${\mathbf P}_x$ be the law of a branching Brownian motion on ${\mathbb R}^d$ 
such that a single particle starts from a point $x\in {\mathbb R}^d$ at time $0$. 
For $t\geq 0$, $x\in {\mathbb R}^d$ and $y\in {\mathbb R}$, 
let $u(t,x,y):={\mathbf P}_x(L_t\leq y)$. 
Then for a fixed $y\in {\mathbb R}$, 
$u(t,x,y)$ as a function of $(t,x)$
is a weak solution to the equation
\begin{equation}\label{eq-pde}
\frac{\partial u}{\partial t}=\frac{1}{2}\Delta u+\left(\sum_{n=0}^{\infty}p_nu^n-u\right)\mu
\quad \text{on $(0,\infty)\times {\mathbb R}^d$}
\end{equation}
with the initial condition 
\begin{equation}\label{eq-initial}
u(0,x,y)={\bf 1}_{\{y\geq |x|\}} \quad \text{for $x\in {\mathbb R}^d$}
\end{equation}
(see \cite[Subsection 1.3]{INW68-1} and \cite[Example 3.4]{INW68-2}). 
Here we recall that $\mu$ and $\{p_n(x)\}_{n\geq 1}$ are  
the branching rate measure and branching mechanism, respectively. 
Corollary \ref{cor-main} and Remark \ref{rem-mtg} imply that 
for $d=1$ or $2$, if $\mu({\mathbb R}^d)<\infty$, then for each $x\in {\mathbb R}^d$,
\begin{equation}\label{eq-lim-1}
\lim_{t\rightarrow\infty}u(t,x,\delta t)=
\begin{cases}
0 & \text{if $0\leq \delta<\sqrt{-\lambda/2}$},\\
1 & \text{if $\delta>\sqrt{-\lambda/2}$}. 
\end{cases}
\end{equation}
Moreover, if $\delta>\sqrt{-\lambda/2}$, 
then it follows by Lemma \ref{lem-positive}  that 
\begin{equation}\label{eq-lim-2}
u(t,x,\delta t)=1-o(e^{-c t}) \quad (t\rightarrow\infty)
\end{equation}
for any $c\in (0, (\lambda+\sqrt{-2\lambda}\delta)\wedge (\delta^2/2))$. 
For $d=1$, a similar assertion is valid for the function $v(t,x,y):={\mathbf P}_x(R_t\leq y)$, 
where the initial condition \eqref{eq-initial} is replaced by 
$$
u(0,x,y)={\bf 1}_{\{y\geq x\}} \quad \text{for $x\in {\mathbb R}$}. 
$$
Related to  the so-called Fisher-Kolmogorov-Petrovskii-Piskunov equation, 
equations like \eqref{eq-lim-1} and \eqref{eq-lim-2} are studied more in detail 
for a class of branching Brownian motions on ${\mathbb R}$  
(see, e.g., \cite{Mc75,Mc76}, \cite{B78}, 
\cite{LS88}, \cite{LS89}, \cite{R13}, \cite{BH16}, 
\cite[Subsection 1.1]{Sh15} and references therein).
\medskip

The rest of this paper is organized as follows. 
In Section \ref{sect-pre-result}, we first introduce the Kato class measures 
and Feynman-Kac semigroups. 
We next introduce the model of branching Brownian motions. 
We finally present our main results and their applications to several concrete models. 
In Section \ref{sect-prf}, we first show results on the growth rate of Feynman-Kac semigroups. 
A few parts of elementary calculations are given in Appendix \ref{appendix-1}. 
We next provide preliminary lemmas for the proofs of Theorems \ref{thm-main} and \ref{thm-direction}. 
We finally complete the proof by showing Lemmas \ref{lem-limsup} and \ref{lem-liminf}.

Throughout this paper, the letters $c$ and $C$ (with subscript) 
denote finite positive constants which may vary from place to place. 
For positive functions $f(t)$ and $g(t)$ on $(0,\infty)$, 
we write $f(t)\sim g(t) \ (t\rightarrow\infty)$ if $f(t)/g(t)\rightarrow 1 \ (t\rightarrow \infty)$. 
We also write $f(t)\lesssim  g(t) \ (t\rightarrow \infty)$ if there exist $c>0$ and $T>0$ such that 
$f(t)\leq cg(t)$ for all $t\geq T$. 
If $f(t)\lesssim g(t) \ (t\rightarrow \infty)$ and $g(t)\lesssim f(t) \ (t\rightarrow \infty)$, 
then we write $f(t)\asymp g(t) \ (t\rightarrow\infty)$.

\section{Preliminaries and results}\label{sect-pre-result}

\subsection{Kato class measures and Feynman-Kac semigroups}
Let ${\bf M}=(\Omega, {\cal F}, \{{\cal F}_t\}_{t\geq 0}, \{B_t\}_{t\geq 0}, \{P_x\}_{x\in {\mathbb R}^d}, 
\{\theta_t\}_{t\geq 0})$ 
be the Brownian motion on ${\mathbb R}^d$, 
where $\{{\cal F}_t\}_{t\geq 0}$ is the minimal admissible filtration 
and $\{\theta_t\}_{t\geq 0}$ is the time shift operator of the paths 
such that $B_s\circ\theta_t=B_{s+t}$ identically for $s,t\geq 0$. 
Let 
$$p_t(x,y)=\frac{1}{(2\pi t)^{d/2}}\exp\left(-\frac{|x-y|^2}{2t}\right) \quad 
(x,y\in {\mathbb R}^d, t>0).$$
Then $p_t(x,y)$ is  the density of the transition function of ${\bf M}$, that is, 
$$P_x(B_t\in A)=\int_A p_t(x,y)\,{\rm d}y, \quad A\in {\cal B}({\mathbb R}^d).$$
For $\alpha>0$, 
let $G_{\alpha}(x,y)$ be the $\alpha$-resolvent density of ${\bf M}$ 
given by 
$$G_{\alpha}(x,y)=\int_0^{\infty}e^{-\alpha t}p_t(x,y)\,{\rm d}t
=\frac{1}{\pi}\left(\frac{\sqrt{2\alpha}}{2\pi |x-y|}\right)^{(d-2)/2}K_{(d-2)/2}(\sqrt{2\alpha}|x-y|),$$
where $K_{\nu}(x)$ is a Bessel function of the third kind:
$$K_{\nu}(t):=\frac{\pi}{2\sin (\nu\pi)}\left(\sum_{n=0}^{\infty}\frac{(t/2)^{-\nu+2n}}{n!\Gamma(-\nu+n+1)}
-\sum_{n=0}^{\infty}\frac{(t/2)^{\nu+2n}}{n!\Gamma(\nu+n+1)}\right)$$
(see, e.g., \cite[Eq.\ 10.25.2 and Eq.\ 10.27.4]{N-DLMF17}). 
For $d=1$ especially,
$$
G_{\alpha}(x,y)=\frac{1}{\sqrt{2\alpha}}e^{-\sqrt{2\alpha} |x-y|}.
$$
For each $\alpha>0$, we also have as $|x-y|\rightarrow 0$,
$$
G_{\alpha}(x,y)\sim 
\begin{cases}
\displaystyle \frac{1}{\sqrt{2\alpha}} & \text{if $d=1$,}\vspace{2mm}\\
\displaystyle \frac{1}{\pi}\log\left(\frac{1}{|x-y|}\right) &\text{if $d=2$,} \vspace{2mm} \\
\displaystyle \frac{\Gamma(d/2-1)}{2\pi^{d/2}}\frac{1}{|x-y|^{d-2}}  & \text{if $d\geq 3$}.
\end{cases} 
$$
(see, e.g., \cite[Example 7.14]{SP14} and \cite[Eq.\ 10.30.2 and Eq.\ 10.30.3]{N-DLMF17}).  	
Noting that 
$$K_{\nu}(t)\sim \sqrt{\frac{\pi}{2t}}e^{-t} \quad (t\rightarrow\infty)$$
(see, e.g., \cite[Eq.\ 10.25.3]{N-DLMF17}), we get for $d\geq 1$,
\begin{equation}\label{eq-resolvent}
G_{\alpha}(x,y)\sim 
\frac{1}{\sqrt{2\alpha}}\left(\frac{\sqrt{2 \alpha}}{2\pi |x-y|}\right)^{(d-1)/2}e^{-\sqrt{2\alpha}|x-y|}
\quad (|x-y|\rightarrow\infty).
\end{equation}
If $d\geq 3$, then we denote by $G(x,y)$ the Green function of ${\bf M}$:
$$G(x,y)=\int_0^{\infty}p_t(x,y)\,{\rm d}t=\frac{\Gamma(d/2-1)}{2\pi^{d/2}}\frac{1}{|x-y|^{d-2}}.$$
We also define $G_0(x,y):=G(x,y)$. 

According to \cite{ABM91, C02, T08},
we first introduce two classes of measures:

\begin{defn}\label{def-kato}
\begin{enumerate}
\item[{\rm (1)}] Let $\mu$ be a positive Radon measure on ${\mathbb R}^d$. 
Then $\mu$ belongs to the Kato class {\rm (}$\mu\in {\cal K}$ in notation{\rm )} 
if one of the following holds:
\begin{enumerate}
\item[{\rm (i)}] $d=1$ and 
$$\sup_{x\in {\mathbb R}}\int_{|x-y|\leq 1}\mu({\rm d}y)<\infty{\rm ;}$$
\item[{\rm (ii)}] $d=2$ and 
$$\lim_{R\rightarrow +0}\sup_{x\in {\mathbb R}^2}\int_{|x-y|\leq R}
\log\left(\frac{1}{|x-y|}\right)\,\mu({\rm d}y)=0{\rm ;}$$
\item[\rm (iii)] $d\geq 3$ and 
$$\lim_{R\rightarrow +0}\sup_{x\in {\mathbb R}^d}\int_{|x-y|\leq R}
G(x,y)\mu({\rm d}y)=0.$$
\end{enumerate}

\item[{\rm (2)}] For $\beta>0$, $\mu\in {\cal K}$ is $\beta$-Green tight 
{\rm (}$\mu\in {\cal K}_{\infty}(\beta)$ in notation{\rm)} if 
$$\lim_{R\rightarrow\infty}\sup_{x\in {\mathbb R}^d}\int_{|y|\geq R}G_{\beta}(x,y)\,\mu({\rm d}y)=0.$$
When $d\geq 3$, $\mu\in {\cal K}$ is Green tight if the equality above is valid 
for $\beta=0$.  
\end{enumerate}
\end{defn}

The Kato class condition imposes  
the uniformly local integrability of measures weighted by the resolvent.
For equivalent definitions of the Kato class, 
see, e.g., \cite[Theorem 3.9 and Remark 3.10]{ABM91}, \cite[Theorem 4.1]{AM92}, 
\cite{KT06, KT07} and references therein.

The $\beta$-Green tight condition restricts the size of measures at infinity 
in terms of the $\beta$-resolvent.
The class ${\cal K}_{\infty}(0)$ was first introduced by Zhao \cite{Z92} 
for transient Brownian motions. 
We see from \cite{T08} that 
${\cal K}_{\infty}(\beta) \ (\beta>0)$ is independent of $\beta$ by the resolvent equation 
and that ${\cal K}_{\infty}(0)\subset {\cal K}_{\infty}(1)$. 
Define
$$
{\cal K}_{\infty}=
\begin{cases}{\cal K}_{\infty}(1) & (d=1,2), \\
{\cal K}_{\infty}(0) & (d\geq 3).
\end{cases}$$
If $\mu$ is a Kato class measure with compact support in ${\mathbb R}^d$, 
then $\mu\in {\cal K}_{\infty}$ by definition. 

\begin{exam}\label{ex-v}\rm  
Let $V(x)$ be a nonnegative Borel function on ${\mathbb R}^d$ 
such that $\mu({\rm d}x):=V(x)\,{\rm d}x$ is a positive Radon measure 
on ${\mathbb R}^d$.
\begin{enumerate}
\item[(i)] Assume that there exist $R_1>1$, $c_1>0$ and $p>0$ 
such that 
$$V(x)\leq \frac{c_1}{|x|^p} \quad \text{for any $x\in {\mathbb R}^d$ with $|x|\leq R_1$}.$$
Then  $\mu\in {\cal K}$ 
if $d=1$ and $p<1$, or if $d\geq 2$ and $p<2$. 
In fact, if $d=1$ and $p<1$, then  
\begin{equation*}
\begin{split}
\int_{|x-y|\leq 1}V(y)\,{\rm d}y
\leq \int_{|x-y|\leq 1}\frac{{\rm d}y}{|y|^p}
&= \int_{|x-y|\leq 1, |x-y|\leq |y|}\frac{{\rm d}y}{|y|^p}
+\int_{|x-y|\leq 1, |x-y|>|y|}\frac{{\rm d}y}{|y|^p}\\
&\leq  \int_{|x-y|\leq 1}\frac{{\rm d}y}{|x-y|^p}
+\int_{|y|\leq 1}\frac{{\rm d}y}{|y|^p}=c
\end{split}
\end{equation*}
for some positive constant $c$ independently of $x$.
A similar calculation applies to the case  $d=2$ and $p<2$.

\item[(ii)] Assume that $\mu\in {\cal K}$. 
If there exist positive constants $R_2$, $c_2$ and $q$ such that 
$$V(x)\leq \frac{c_2}{|x|^q} \quad \text{for any $x\in {\mathbb R}^d$ with $|x|\geq R_2$},$$
then by the same way as in \cite[Proposition 2]{Z92}, 
we have $\mu\in {\cal K}_{\infty}$ if $d=1$ and $q>1$, or if $d\geq 2$ and $q>2$. 
\end{enumerate}
\end{exam}

\begin{exam}\label{ex-kato-1d}\rm 
For $d=1$, any finite measure $\mu$ belongs to ${\cal K}_{\infty}$ 
because 
$$\sup_{x\in {\mathbb R}}\int_{|x-y|\leq 1}\mu({\rm d}y)\leq \mu({\mathbb R})<\infty$$
and 
$$\sup_{x\in {\mathbb R}}\int_{|y|\geq R}G_1(x,y)\,\mu({\rm d}y)\leq \frac{1}{\sqrt{2}}\mu(|x|\geq R)\rightarrow 0 \quad 
(R\rightarrow \infty).$$
\end{exam}

\begin{exam}\label{ex-kato-surface}\rm 
For $d\geq 2$, let $\delta_R$ be the surface measure on the sphere 
$\{x\in {\mathbb R}^d \mid |x|=R\}$ for $R>0$. 
Then $\delta_R$ belongs to ${\cal K}_{\infty}$ because 
it is a Kato class measure with compact support by \cite[Theorem 4.1]{BEKS94}.
\end{exam}

We next introduce the notion of  positive continuous additive functionals.
Let $A=\{A_t\}_{t\geq 0}$ be a $[0,\infty]$-valued stochastic process on $(\Omega,{\cal F})$. 
We say that $A$ is a {\it positive continuous additive functional {\rm (}in the strict sense{\rm )}}  
(PCAF in short) of ${\mathbf M}$ if 
\begin{enumerate}
\item $A_t$ is ${\cal F}_t$-measurable for any $t\geq 0$;
\item There exists an event $\Lambda\in {\cal F}_{\infty}$, which satisfies 
$P_x(\Lambda)=1$ for any $x\in {\mathbb R}^d$ and 
$\theta_t\Lambda\subset \Lambda$ for any $t>0$, 
such that for any $\omega\in \Lambda$, 
\begin{itemize}
\item $A_0(\omega)=0$;
\item $A_t(\omega)$ is finite and continuous in $t\in [0,\infty)$;
\item $A_{t+s}(\omega)=A_t(\omega)+A_s(\theta_t\omega)$ for any $s,t\geq 0$
\end{itemize}
\end{enumerate}
(see, e.g., \cite{ABM91} and \cite[p.401]{FOT11}).
For each $\mu\in {\cal K}$, there exists 
a unique PCAF ($A^{\mu}$ in notation) 
such that for any nonnegative Borel function $f$, 
$$\lim_{t\rightarrow 0}\frac{1}{t}\int_{{\mathbb R}^d}E_x\left[\int_0^t f(B_s)\,{\rm d}A_s^{\mu}\right]\,{\rm d}x
=\int_{{\mathbb R}^d}f(x)\,\mu({\rm d}x)$$
(\cite[Proposition 3.8]{ABM91} and  \cite[Theorems 5.1.3 and  5.1.7]{FOT11}). 
We note that if $d\geq 3$, 
then by \cite[Proposition 2.2]{C02}, any measure $\mu\in {\cal K}_{\infty}$ is {\it Green-bounded}:
\begin{equation}\label{eq-green-bd}
\sup_{x\in {\mathbb R}^d}
E_x\left[A_{\infty}^{\mu}\right]=\sup_{x\in {\mathbb R}^d} \int_{{\mathbb R}^d}G(x,y)\,\mu({\rm d}y)<\infty.
\end{equation}

If $\mu({\rm d}x)=g(x)\,{\rm d}x$ for a nonnegative Borel function $g$ on ${\mathbb R}^d$,  
then $A_t^{\mu}=\int_0^t g(B_s)\,{\rm d}s$. 
For $d=1$, the Dirac measure at $a\in {\mathbb R}$ ($\delta_a$ in notation)  
belongs to ${\cal K}_{\infty}(1)$ by Example \ref{ex-kato-1d} 
and $A_t^{\delta_a}=2L(t,a)$, where $L(t,a)$ is the Brownian local time at $a$ 
given by  
$$L(t,a)=\lim_{\varepsilon\rightarrow 0}
\frac{1}{4\varepsilon}\int_0^t{\bf 1}_{\{B_s\in (a-\varepsilon,a+\varepsilon)\}}\,{\rm d}s$$
(see, e.g., \cite[Example 5.1.1]{FOT11} and \cite[Chapter III, Section 4]{IW89}).
\medskip

For a fixed measure  $\mu\in {\cal K}$, 
the multiplicative functional $e^{A_t^{\mu}}$ is called the {\it Feynman-Kac functional}. 
Using it, we finally define the {\it Feynman-Kac semigroup} $\{p_t^{\mu}\}_{t\geq 0}$ by 
$$p_t^{\mu}f(x):=E_x\left[e^{A_t^{\mu}}f(B_t)\right], \quad f\in L^2({\mathbb R}^d)\cap {\cal B}_b({\mathbb R}^d),$$
where ${\cal B}_b({\mathbb R}^d)$ stands for the set of bounded Borel functions on ${\mathbb R}^d$.
Then $\{p_t^{\mu}\}_{t\geq 0}$ forms  
a strongly continuous symmetric semigroup on $L^2({\mathbb R}^d)$ 
such that the $L^2$-generator is formally expressed  as 
the Schr\"odinger type operator ${\cal H}^{\mu}:=-\Delta/2-\mu$.  
We can extend $\{p_t^{\mu}\}_{t\geq 0}$ to $L^p({\mathbb R}^d)$ 
for any $p\in [1,\infty]$ (\cite[Theorem 6.1 (i)]{ABM91}). 
For simplicity, we use the same notation for the extended semigroup. 
Let $\|\cdot\|_{p,q}$ denote the operator norm from $L^p({\mathbb R}^d)$ to $L^q({\mathbb R}^d)$. 
We then have

\begin{thm}\label{thm-abm} {\rm (\cite[Theorems 6.1 (iii) and 7.1 (ii)]{ABM91})} Let $\mu\in {\cal K}$.
\begin{enumerate}
\item For any $t>0$, $\|p_t^{\mu}\|_{p,q}<\infty$ for any $1\leq p\leq q\leq \infty$.
\item For any $f\in {\cal B}_b({\mathbb R}^d)$ and $t>0$, $p_t^{\mu}f$ is a bounded continuous function on ${\mathbb R}^d$.
\end{enumerate}
\end{thm}

Assume that $\mu\in {\cal K}_{\infty}(1)$.   
Define 
$$\lambda({\mu}):=\inf\left\{\frac{1}{2}\int_{{\mathbb R}^d}|\nabla u|^2\,{\rm d}x-\int_{{\mathbb R}^d}u^2\,{\rm d}\mu 
\mid u\in C_0^{\infty}({\mathbb R}^d), \int_{{\mathbb R}^d}u^2\,{\rm d}x=1\right\},$$
where $C_0^{\infty}({\mathbb R}^d)$ stands for the set of smooth functions on ${\mathbb R}^d$ with compact support. 
Then $\lambda(\mu)$ is the bottom of the $L^2$-spectrum of ${\cal H}^{\mu}$. 
In particular,  
if $\lambda(\mu)<0$, then $\lambda(\mu)$ is the principal eigenvalue of ${\cal H}^{\mu}$ 
(see \cite[Lemma 4.3]{T03} or \cite[Theorem 2.8]{T08}) 
and the corresponding eigenfunction $h$ 
has a bounded, continuous and strictly positive version by  Theorem \ref{thm-abm}
(see, e.g., \cite[Section 4]{T08}). 

\begin{rem}\label{rem-positive}\rm 
Let $\mu\in {\cal K}_{\infty}(1)$ be non-trivial. 
Then by \cite[Lemmas 4.1 and 4.2]{T03}, 
\begin{enumerate}
\item[(i)] If $d=1$ or $2$, then $\lambda(\mu)<0$. 
\item[(ii)] If $d\geq 3$, then there exists $c_*>0$ such that 
$\lambda(c \mu)<0$ if and only if $c>c_*$. 
\end{enumerate}
\end{rem}

In what follows, we assume that the eigenfunction $h$ is 
bounded, continuous and strictly positive on ${\mathbb R}^d$
such that $\int_{{\mathbb R}^d}h(x)^2\,{\rm d}x=1$. 
Then by the proof of \cite[Theorem 5.2]{T08}, we have for any $f\in {\cal B}_b({\mathbb R}^d)$,
\begin{equation}\label{eq-fkf-asymp}
\lim_{t\rightarrow\infty}e^{\lambda t}E_x\left[e^{A_t^{\mu}}f(B_t)\right]
=h(x)\int_{{\mathbb R}^d}f(y)h(y)\,{\rm d}y, \quad x\in {\mathbb R}^d.
\end{equation}
If $\mu$ is further compactly supported, 
then by \eqref{eq-resolvent} and \cite[Remark 2.4]{S08}, 
\begin{equation}\label{eq-ground}
h(x)\asymp G_{-\lambda}(0,x)\asymp \frac{e^{-\sqrt{-2\lambda}|x|}}{|x|^{(d-1)/2}}  \quad (|x|\rightarrow\infty).
\end{equation}

\subsection{Branching Brownian motions}\label{subsect-model}
In this subsection, we introduce the model of branching Brownian motions 
by following \cite{INW68-1,INW68-2, INW69}. 
For $x\in {\mathbb R}^d$, let $\{p_n(x)\}_{n\geq 1}$ be a sequence such that 
$$0\leq p_n(x)\leq 1  \quad (n\geq 1) \quad \text{and} 
\quad \sum_{n=1}^{\infty}p_n(x)=1.$$
Let $\tau$ be the nonnegative random variable defined on $(\Omega, {\cal F}, P_x)$, 
which is independent of the Brownian motion, of exponential distribution 
with rate $1$; $P(\tau>t)=e^{-t}$ for any $t>0$. 
Let $\mu$ be a Kato class measure on ${\mathbb R}^d$.
We define 
$$Z:=\inf\left\{t>0 \mid A_t^{\mu}\geq \tau\right\}$$
so that 
$$P_x(Z>t \mid {\cal F}_{\infty})=e^{-A_t^{\mu}}.$$

We can describe the branching Brownian motion as follows: 
a Brownian particle $\{B_t\}_{t\geq 0}$ starts at $x\in {\mathbb R}^d$ 
according to the law $P_x$. 
At time $Z$, this particle splits into $n$ particles with probability $p_n(B_Z)$ ($n\geq 1$). 
These particles then start at $B_Z$ independently according to the law $P_{B_Z}$, 
and each of them continues the same procedure. 
Here we note that a particle gives birth to at least one particle at the branching time 
so that non-extinction occurs.

Let us define the equivalent relation $\sim$ 
on $({\mathbb R}^d)^n=\underbrace{{\mathbb R}^d\times \dots\times {\mathbb R}^d}_{n}$ as follows:
for ${\mathbf x}^n=(x^1,\dots, x^n)$ and ${\mathbf y}^n=(y^1,\dots, y^n)\in ({\mathbb R}^d)^n$, 
we write ${\mathbf x}\sim {\mathbf y}$ if there exists a permutation $\sigma$ on $\{1,2,\dots, n\}$ 
such that $y^i=x^{\sigma(i)}$ for any $i\in \{1,2,\dots, n\}$. 
If we define  $({\mathbb R}^d)^{(n)}=({\mathbb R}^d)^n/\sim$ 
and ${\mathbf X}=\cup_{n=1}^{\infty}({\mathbb R}^d)^{(n)}$, 
then $n$ points in ${\mathbb R}^d$ determine a point in $({\mathbb R}^d)^{(n)}$.
Hence we can define 
the branching Brownian motion 
$\overline{\mathbf M}=(\{{\mathbf B}_t\}_{t\geq 0}, \{P_{\mathbf x}\}_{{\mathbf x}\in {\mathbf X}})$ 
on ${\mathbf X}$ with branching rate $\mu$ and branching mechanism $\{p_n(x)\}_{n\geq 1}$.  

Let $T$ be the first splitting time of $\overline{\mathbf M}$ given by 
\begin{equation}\label{eq-split}
{\mathbf P}_x(T>t \mid \sigma(B))=P_x(Z>t \mid {\cal F}_{\infty})=e^{-A_t^{\mu}} \quad (t>0).
\end{equation}
By definition, the first splitting time becomes small if the particle 
moves on the support of $\mu$ often. 
Let 
$$Q(x):=\sum_{n=1}^{\infty}np_n(x) (\geq 1)$$
be the expected offspring number at branching cite $x\in {\mathbb R}^d$. 
Denote by $Z_t$ the total number of particles at time $t$, that is, 
$$Z_t=n \quad 
\text{if ${\mathbf B}_t=({\mathbf B}_t^{1},\dots, {\mathbf B}_t^n)
\in {({\mathbb R}^d)}^{(n)}$.}$$
We define for $f\in {\cal B}_b({\mathbb R}^d)$, 
$$Z_t(f):=\sum_{k=1}^{Z_t}f({\mathbf B}_t^k).$$
For $A\in {\cal B}({\mathbb R}^d)$, 
let $Z_t(A):=Z_t({\bf 1}_A)$ be the number of particles in the set $A$ at time $t$.
If the measure 
$$
\nu({\rm d}x):=(Q(x)-1)\mu({\rm d}x)
$$
belongs to the Kato class, 
then by the same way as in of \cite[Lemma 3.3]{S08}, 
we have 
\begin{equation}\label{eq-br-fk}
{\mathbf E}_x\left[Z_t(f)\right]
=E_x\left[e^{A_t^{\nu}}f(B_t)\right].
\end{equation}

Assume that $\nu\in {\cal K}_{\infty}(1)$ and $\lambda:=\lambda(\nu)<0$. 
Let $h$ be the eigenfunction of ${\cal H}^{(Q-1)\mu}$ corresponding to $\lambda$ and 
$$M_t:=e^{\lambda t}Z_t(h) \quad (t\geq 0).$$
Since $M_t$ is a nonnegative ${\mathbf P}_x$-martingale,
the limit $M_{\infty}:=\lim_{t\rightarrow\infty}M_t\in [0,\infty)$ exists 
${\mathbf P}_x$-a.s. 
Furthermore, by \cite[Theorem 3.7]{CS07}, there exists an event of ${\mathbf P}_x$-full probability measure 
such that under this event,   
\begin{equation}\label{eq-limit-thm}
\lim_{t\rightarrow\infty}e^{\lambda t}Z_t(A)=M_{\infty}\int_A h(y)\,{\rm d}y
\end{equation}
for any $A\in {\cal B}({\mathbb R}^d)$ 
such that its boundary has zero Lebesgue measure.

\begin{rem}\label{rem-llogl}\rm 
Suppose that 
\begin{equation}\label{eq-llogl}
\sup_{x\in {\mathbb R}^d}\sum_{n=2}^{\infty}(n\log n) p_n(x)<\infty.
\end{equation}
Then by \cite[Theorem 1.2]{CRY17+}, 
$M_t$ converges to $M_{\infty}$ in $L^1({\mathbf P}_x)$.  
This implies that  ${\mathbf E}_x[M_{\infty}]=h(x)$ and thus ${\mathbf P}_x(M_{\infty}>0)>0$. 
\begin{itemize}
\item For $d=1$ or $2$,  
if we further assume that $\mu({\mathbb R}^d)<\infty$, 
then ${\mathbf P}_x(M_{\infty}>0)=1$ by the same way as 
in \cite[Proposition 3.6 and  Remark 3.14]{S08}. 
\item For $d\geq 3$, if we further assume that $\mu\in {\cal K}_{\infty}$, 
then ${\mathbf P}_x(M_{\infty}=0)>0$. In fact, 
since $\nu\in {\cal K}_{\infty}$, \eqref{eq-green-bd} yields that  
$$\sup_{x\in {\mathbb R}^d}E_x[A_{\infty}^{\nu}]<\infty,$$
that is, $P_x(A_{\infty}^{\nu}<\infty)=1$.  
Moreover, noting that  $\{M_{\infty}=0\}\supset \{T=\infty\}$, we get
$${\mathbf P}_x(M_{\infty}=0)\geq {\mathbf P}_x(T=\infty)=E_x[e^{-A_{\infty}^{\nu}}]>0.$$
\end{itemize}
\end{rem}

\subsection{Results}
For $R>0$, let $B_R:=\{x\in {\mathbb R}^d \mid |x|<R\}$ be an open ball 
with radius $R$ centered at the origin. 
Denote by $Z_t^R$ the number of particles outside $B_R^c$ at time $t$,   
that is, $Z_t^R=Z_t(B_R^c)$.  
We then obtain

\begin{thm}\label{thm-main}
Assume that the measure $\nu_{\beta}({\rm d}x):=e^{\beta|x|}\nu({\rm d}x)$ belongs to 
${\cal K}_{\infty}(1)$ for any $\beta>0$. If $\lambda<0$, 
then the following assertions hold:
\begin{enumerate}
\item[{\rm (i)}] If $\delta>\sqrt{-\lambda/2}$, then 
$$\lim_{t\rightarrow\infty} Z_t^{\delta t}=0,
\quad \text{${\mathbf P}_x$-a.s.}$$ 
\item[{\rm (ii)}] If ${\mathbf P}_x(M_{\infty}>0)>0$, then for any $\delta\in [0,\sqrt{-\lambda/2})$,
$$
\lim_{t\rightarrow\infty}\frac{1}{t}\log Z_t^{\delta t}=-\lambda-\sqrt{-2\lambda}\delta,
\quad \text{${\mathbf P}_x(\cdot\mid M_{\infty}>0)$-a.s.}$$
\end{enumerate}
\end{thm}

Theorem \ref{thm-main} says that if $\delta>\sqrt{-\lambda/2}$, 
then all the particles at time $t$ will be inside the ball $B_{\delta t}$ eventually. 
On the other hand, if $\delta<\sqrt{-\lambda/2}$, 
then for all sufficiently large $t>0$, 
there will be particles outside the ball $B_{\delta t}$ at time $t$ 
such that the number of such particles grows exponentially.  
Note that Theorem \ref{thm-main} (ii) is valid for $\delta=0$ without the condition on $\nu_{\beta}$ 
(Lemma \ref{lem-growth}). 

Under the assumption in Theorem \ref{thm-main}, 
$\lambda$ is the principal eigenvalue 
of the Schr\"odinger type operator ${\cal H}^{(Q-1)\mu}$
and thus the associated eigenfunction $h$ exists. 
We impose the technical condition on $\nu_{\beta}$ 
in order to verify the growth rate of the expectation of $Z_t^{\delta t}$ 
(see Subsection \ref{subsect-fk} below). 
For instance, this condition is satisfied if $\nu$ is a Kato class measure 
with compact support in ${\mathbb R}^d$.  
We will prove Theorem \ref{thm-main} in Section \ref{sect-prf}. 
\medskip

By Theorem \ref{thm-main} (ii), 
we can determine the spread rate of branching Brownian motions.  
Let 
$$L_t:=\max_{1\leq k\leq Z_t}|{\mathbf B}_t^k|$$ 
be the maximum of the Euclidean norms of particles alive at time $t$. 
We then have 
\begin{cor}\label{cor-main}
Under the same conditions as in Theorem {\rm \ref{thm-main}}, 
if ${\mathbf P}_x(M_{\infty}>0)>0$, then 
\begin{equation}\label{eq-cor}
\lim_{t\rightarrow\infty}\frac{L_t}{t}=\sqrt{-\frac{\lambda}{2}}, 
\quad \text{${\mathbf P}_x(\cdot\mid M_{\infty}>0)$-a.s.}
\end{equation}
\end{cor}

\pf \ We first let $\delta>\sqrt{-\lambda/2}$.  
Since 
$${\mathbf P}_x\left(\text{$L_t\leq \delta t$ for all sufficiently large $t$}\right)=1$$
by Theorem \ref{thm-main} (i), we obtain   
$$\limsup_{t\rightarrow\infty}\frac{L_t}{t}\leq \delta, \quad 
\text{${\mathbf P}_x$-a.s.}$$
By letting $\delta\searrow\sqrt{-\lambda/2}$, we get 
$$\limsup_{t\rightarrow\infty}\frac{L_t}{t}\leq \sqrt{-\frac{\lambda}{2}}, \quad \text{${\mathbf P}_x$-a.s.}$$

We next let $\delta\in (0,\sqrt{-\lambda/2})$.  In what follows, 
we assume that ${\mathbf P}_x(M_{\infty}>0)>0$. Since
$${\mathbf P}_x\left(\lim_{t\rightarrow\infty}Z_t^{\delta t}=\infty\mid M_{\infty}>0\right)=1$$
by Theorem \ref{thm-main} (ii), we have
$${\mathbf P}_x\left(\text{$L_t\geq \delta t$ for all sufficiently large $t$}\mid M_{\infty}>0\right)=1,$$
which yields that 
$$\liminf_{t\rightarrow\infty}\frac{L_t}{t}\geq \delta, \quad \text{${\mathbf P}_x(\cdot\mid M_{\infty}>0)$-a.s.}$$
By letting $\delta\nearrow\sqrt{-\lambda/2}$, we get 
$$\liminf_{t\rightarrow\infty}\frac{L_t}{t}\geq \sqrt{-\frac{\lambda}{2}}, \quad \text{${\mathbf P}_x(\cdot\mid M_{\infty}>0)$-a.s.}$$
Hence the proof is complete. 
\qed
\medskip

By the same way as for the proofs of Theorem \ref{thm-main} and Corollary \ref{cor-main}, 
we can also show that the spread rate is uniform in direction.
For a unit vector $r$ in ${\mathbb R}^d$, define 
$B_R^r:=\{x\in {\mathbb R}^d \mid \langle x, r\rangle <R\}$ for $R>0$ and 
$Z_t^{\delta t,r}:=Z_t((B_{\delta t}^r)^c)$. 
Let 
$$L_t^r:=\max_{1\leq k\leq Z_t}\langle {\mathbf B}_t^k,r\rangle$$ 
be 
the maximal displacement in direction $r$ of particles alive at time $t$,  
and let $K_r(t)$ be the index of a particle at time $t$ 
such that $L_t^r=\langle {\mathbf B}_t^{K_r(t)},r\rangle$.
 
\begin{thm}\label{thm-direction}\rm 
Let $r$ be a unit vector in ${\mathbb R}^d$. 
Assume  the same conditions as in Theorem {\rm \ref{thm-main}}. 
\begin{enumerate}
\item[(i)] The assertions of Theorem \ref{thm-main} are valid 
by replacing $Z_t^{\delta t}$ with $Z_t^{\delta t,r}$. 
\item[(ii)] If ${\mathbf P}_x(M_{\infty}>0)>0$, then 
\begin{equation}\label{eq-direction-0}
\lim_{t\rightarrow\infty}\frac{L_t^r}{t}=\sqrt{-\frac{\lambda}{2}}, 
\quad \text{${\mathbf P}_x(\cdot \mid M_{\infty}>0)$-a.s.}
\end{equation}
and
\begin{equation}\label{eq-direction}
\lim_{t\rightarrow\infty}\frac{{\mathbf B}_t^{K_r(t)}}{t}=\sqrt{-\frac{\lambda}{2}}r, 
\quad \text{${\mathbf P}_x(\cdot \mid M_{\infty}>0)$-a.s.}
\end{equation}
\end{enumerate}
\end{thm}

We will also prove Theorem \ref{thm-direction} in Section \ref{sect-prf}.

\begin{rem}\label{rem-mtg}\rm  
Assume that $d=1$ or $2$. 
If $\mu({\mathbb R}^d)<\infty$ and \eqref{eq-llogl} holds, 
then ${\mathbf P}_x(M_{\infty}>0)=1$ by Remark \ref{rem-llogl}. 
Hence if we further assume that $\nu_{\beta}\in {\cal K}_{\infty}(1)$ for any $\beta>0$, 
then by Remark \ref{rem-positive}, 
the assertions of Theorem \ref{thm-main}, Corollary \ref{cor-main} and Theorem \ref{thm-direction} 
are valid ${\mathbf P}_x$-a.s.\ in place of ${\mathbf P}_x(\cdot\mid M_{\infty}>0)$-a.s. 
Note that if $d=1$ and $r=1$, then $L_t^1=\sup_{1\leq k\leq Z_t}{\mathbf B}_t^k$.
\end{rem}

We also discuss the spread rate of $\overline{\bf M}$
under the condition that $\lambda=0$. 
We first note that if $\lambda=0$, then $d\geq 3$ by Remark \ref{rem-positive}, 
that is, ${\mathbf M}$ is transient. 
Furthermore, if we also assume that 
$\sup_{x\in {\mathbb R}^d}Q(x)<\infty$, then  $\overline{\mathbf M}$ exhibits {\it local extinction}: 
for any relatively compact open set $A\subset {\mathbb R}^d$, 
$${\mathbf P}_x\left(\text{there exists $S>0$ such that $Z_t(A)=0$ for all $t\geq S$}\right)=1$$
(see, e.g., \cite[Theorem 3.2]{S06} and \cite[Theorem 3.8]{S08}). 
Hence if the measure $\nu$ is compactly supported in ${\mathbb R}^d$, 
then for all sufficiently large time, 
all the particles will be outside the support of $\nu$, 
that is, branching will never occur.
As a result, $L_t$ becomes the maximal Euclidean norm 
among the finite number of Brownian particles so that  
it satisfies the law of the iterated logarithm. 
In order to show this,  
we recall that $T$ is the first splitting time of $\overline{\bf M}$. 
Let $\theta_t$ be the time shift operator of the paths for $\overline{\bf M}$ and 
$$T_1:=T, \quad T_n:=T_{n-1}+T\circ\theta_{T_{n-1}}  \quad (n\geq 2).$$ 
Then $T_n$ denotes the $n$th splitting time of $\overline{\bf M}$.

\begin{prop}\label{prop-zero}
Assume that $\sup_{x\in {\mathbb R}^d}Q(x)<\infty$. 
If $\mu$ is a Kato class measure with compact support in ${\mathbb R}^d$ 
and $\lambda=0$, then 
\begin{equation}\label{eq-lil}
\limsup_{t\rightarrow\infty}\frac{L_t}{\sqrt{2t\log\log t}}=1, 
\quad \text{${\mathbf P}_x$-a.s.}
\end{equation}
\end{prop}

\pf \  Under the full condition of this proposition, 
since $\overline{\mathbf M}$ exhibits local extinction and 
$\mu$ is compactly supported, 
the total number of branching for $\overline{\bf M}$ is finite 
so that 
\begin{equation}\label{eq-finite}
{\mathbf P}_x
\left(\{T=\infty\}\cup \bigcup_{n=1}^{\infty} \left\{T_n<\infty, T_{n+1}=\infty \right\}\right)=1.
\end{equation}
Hence if we define the event $A$ by 
$$A:=\left\{\limsup_{t\rightarrow\infty}\frac{L_t}{\sqrt{2t\log\log t}}=1\right\},$$  
then 
\begin{equation}\label{eq-a-1}
\begin{split}
{\mathbf P}_x(A)
&={\mathbf P}_x\left(\left(\{T=\infty\}\cup \bigcup_{n=1}^{\infty}\left\{T_n<\infty, T_{n+1}=\infty\right\}\right)\cap  A\right)\\
&={\mathbf P}_x\left(\{T=\infty\}\cap A\right)
+\sum_{n=1}^{\infty}{\mathbf P}_x\left(\{T_n<\infty, T_{n+1}=\infty\}\cap  A\right).
\end{split}
\end{equation}

By Khintchine's law of the iterated logarithm, 
\begin{equation}\label{eq-a-inf}
{\mathbf P}_x\left(\{T=\infty\}\cap A\right)
=E_x\left[e^{-A_{\infty}^{\mu}};
\limsup_{t\rightarrow\infty}\frac{|B_t|}{\sqrt{2t\log\log t}}= 1\right]
=E_x\left[e^{-A_{\infty}^{\mu}}\right]={\mathbf P}_x(T=\infty).
\end{equation}
The last probability above is positive by Remark \ref{rem-llogl}. 
For $n\geq 1$, we have by the strong Markov property of $\overline{\mathbf M}$,  
\begin{equation}\label{eq-a-smp}
{\mathbf P}_x\left(\{T_n<\infty, T_{n+1}=\infty\}\cap  A\right)
={\mathbf E}_x\left[{\mathbf P}_{{\mathbf B}_{T_n}}\left(\{T=\infty\}\cap A\right);T_n<\infty\right].
\end{equation}
Since 
$${\mathbf P}_{{\mathbf x}^n}\left(\{T=\infty\}\cap A\right)
\geq \prod_{k=1}^n{\mathbf P}_{x^k}\left(\{T=\infty\}\cap A\right)
=\prod_{k=1}^n{\mathbf P}_{x^k}\left(T=\infty \right)$$
by \eqref{eq-a-inf}, the right hand side of \eqref{eq-a-smp} is greater than  
$${\mathbf E}_x\left[\prod_{k=1}^{Z_{T_n}}
{\mathbf P}_{{\mathbf B}_{T_n}^k}\left(T=\infty\right);T_n<\infty\right]\\
={\mathbf P}_x\left(T_n<\infty, T_{n+1}=\infty\right).
$$
Then by \eqref{eq-finite}, the last expression of \eqref{eq-a-1} is not less than $1$,
which completes the proof. 
\qed

\subsection{Examples}\label{subsect-exam}
In this subsection, we apply the results in the previous subsection 
to several concrete models.  
Throughout this section, we assume that $p_2(x)\equiv 1$  for simplicity. 
Then $\nu=\mu$. 

\begin{exam}\rm 
Assume that $d=1$. 
Let $\mu$ be any branching rate measure such that for any $\beta>0$, 
$\nu_{\beta}({\rm d}x):=e^{\beta |x|}\mu({\rm d}x)$ is a finite measure.  
Then $\nu_{\beta}\in {\cal K}_{\infty}(1)$ and $\lambda<0$ 
by Example \ref{ex-kato-1d} and Remark \ref{rem-positive},  respectively. 
Therefore, the assertions of Theorem \ref{thm-main}, 
Corollary \ref{cor-main} and Theorem \ref{thm-direction} are valid 
${\mathbf P}_x$-a.s.\ by Remark \ref{rem-mtg}. 

\begin{enumerate}
\item Let $\mu=c\delta_0$ for some $c>0$. 
Since $\lambda=-c^2/2$ (see, e.g., \cite[Example 4.4]{S08}), 
we see that if $\delta\in (0,c/2)$, then 
$$\lim_{t\rightarrow\infty}\frac{1}{t}\log Z_t^{\delta t}=\frac{c^2}{2}-c\delta,
\quad \text{${\mathbf P}_x$-a.s.}$$
We also have 
$$\lim_{t\rightarrow\infty}\frac{L_t}{t}=\lim_{t\rightarrow\infty}\frac{R_t}{t}=\frac{c}{2}, \quad 
\text{${\mathbf P}_x$-a.s.}$$
These results were proved by \cite{BH14}. 
\item Suppose that $\mu=c_1\delta_{-a}+c_2\delta_a$ 
for some positive constants $c_1$ and $c_2$ with $c_1\leq c_2$. 
Then the assertions of Theorem \ref{thm-main}, 
Corollary \ref{cor-main} and Theorem \ref{thm-direction} are valid 
${\mathbf P}_x$-a.s. 
Here $\lambda$ is a unique solution to 
$$(c_1-\sqrt{-2\lambda})(c_2-\sqrt{-2\lambda})=c_1c_2e^{-4a\sqrt{-2\lambda}}$$
under the restriction that $\lambda<-(c_2)^2/2$ 
(see, e.g., \cite[p.144]{AGHH88}).

\item 
Suppose that 
$$\mu({\rm d}x)
=\sum_{n=-\infty}^{\infty}e^{-|x|^p}{\delta_{n}}({\rm d}x)
=\sum_{n=-{\infty}}^{\infty}e^{-|n|^p}\delta_n({\rm d}x)$$
for some $p>1$. 
Then the assertions of Theorem \ref{thm-main}, Corollary \ref{cor-main} and Theorem \ref{thm-direction} are valid 
${\mathbf P}_x$-a.s.
By \cite[(2.6)]{S08}, the principal eigenvalue $\lambda$ and corresponding eigenfunction $h$ are characterized by 
the equation
$$h(x)=\frac{1}{\sqrt{-2\lambda}}\sum_{n=-\infty}^{\infty}h(n)e^{-|n|^p-\sqrt{-2\lambda}|x-n|} 
\left(=\int_{{\mathbb R}}G_{-\lambda}(x,y)h(y)\,\mu({\rm d}y)\right) 
\quad (x\in {\mathbb R}).$$
\end{enumerate}
\end{exam}

\begin{exam}\rm Assume that $d\geq 2$. 
For $R>0$, let $\delta_R$ be the surface measure 
on the sphere $\{x\in {\mathbb R}^d \mid |x|=R\}$. 
For $c>0$, we take $\mu:=c\delta_R$ as the branching rate measure. 
Since $\mu$ is a Kato class measure with compact support 
as mentioned in Example \ref{ex-kato-surface}, 
we can apply the results in the previous subsection to our model. 

If $d=2$, then for any $c>0$ and $R>0$, 
the assertions of Theorem \ref{thm-main}, Corollary \ref{cor-main} and Theorem \ref{thm-direction} are valid 
${\mathbf P}_x$-a.s.\ by Remark \ref{rem-mtg}. 
On the other hand, if $d\geq 3$, 
then by \cite[Lemma 3.5]{T02} and \cite[Exercise 6.4.13]{FOT11},
$$\lambda<0 \iff cR>\frac{d-2}{2}.$$
Hence we have \eqref{eq-cor} and \eqref{eq-direction} for $cR>(d-2)/2$ and \eqref{eq-lil} otherwise. 
For $d=3$, we know from \cite[(2.39)]{AGS87} that if $cR>1/2$, 
then $\lambda$ is a unique solution to the equation 
$$\frac{2\sqrt{-2\lambda}e^{2\sqrt{-2\lambda}R}}{e^{2\sqrt{-2\lambda}R}-1}=c.$$
\end{exam}

\begin{exam}\rm 
Let $V(x)$ be a nonnegative function on ${\mathbb R}^d$. 
For $c>0$, we take $\mu({\rm d}x):=cV(x)$ as the branching rate measure.

\begin{enumerate}
\item[(i)] Let $V(x)={\bf 1}_{\{|x|\leq R\}}(x)$ for some $R>0$.
Since $\mu$ is a Kato class measure with compact support by Example \ref{ex-v}, 
we can apply the results in the previous subsection to our model. 

If $d=1$ or $2$, then for any $c>0$ and $R>0$, 
\eqref{eq-cor} and \eqref{eq-direction} hold ${\mathbf P}_x$-a.s.\ by Remark \ref{rem-mtg}. 
On the other hand, if $d\geq 3$, then Remark \ref{rem-positive} implies  
that for each fixed $R>0$, 
there exists $c_*>0$ such that $\lambda<0$ if and only if $c>c_*$. 
Hence we get \eqref{eq-cor} and \eqref{eq-direction} for $c>c_*$ and \eqref{eq-lil} otherwise. 
For $d=3$, we  get  $c_*=\pi^2/(8R^2)$ by a direct calculation.  
Furthermore, for each $c>c_*$, $\lambda\in (c_*-c, 0)$ is a unique solution to the equation
$$\frac{\tan(\sqrt{2(\lambda+c)}R)}{\sqrt{2(\lambda+c)}R}=-\frac{1}{\sqrt{-2\lambda}}.$$

\item[(ii)]
Let $V(x)$ be a nonnegative function with compact support 
such that $V\not\equiv 0$ and $V(x)\leq {\bf 1}_{\{0<|x|\leq R\}}|x|^{-p}$ 
for some $R>0$ and $p\in {\mathbb R}$. 
If $d=1$ and $p<1$, or if $d\geq 2$ and $p<2$, 
then $\mu$ is a Kato class measure with compact support by Example \ref{ex-v}. 
Accordingly, we can apply the results in the previous subsection to our model. 

If $d=1$ and $p<1$, or if  $d=2$ and $p<2$, then 
\eqref{eq-cor} and \eqref{eq-direction} hold ${\mathbf P}_x$-a.s.\ by Remark \ref{rem-mtg}. 
On the other hand, if $d\geq 3$ and $p<2$, 
then Remark \ref{rem-positive} yields that for each fixed $R>0$, 
there exists $c_*>0$ such that $\lambda<0$ if and only if $c>c_*$.
Therefore,  we obtain \eqref{eq-cor} and \eqref{eq-direction} for $c>c_*$ and \eqref{eq-lil} otherwise.

\item[(iii)] 
Let $V(x)$ be a positive function such that $V(x)\leq e^{-|x|^p}$ for some $p>1$. 
Then for any $\beta>0$, 
$\nu_{\beta}({\rm d }x):=e^{\beta|x|}\mu({\rm d}x)$ belongs ${\cal K}_{\infty}(1)$ by Example \ref{ex-v}. 
We can thus apply Theorem \ref{thm-main}, Corollary \ref{cor-main} and Theorem \ref{thm-direction} 
to our model. 

If $d=1$ or $2$, then \eqref{eq-cor} and \eqref{eq-direction} hold ${\mathbf P}_x$-a.s.\ by Remark \ref{rem-mtg}. 
On the other hand, if $d\geq 3$, then Remark \ref{rem-positive} implies that there exists $c_*>0$ such that 
$\lambda<0$ if and only if $c>c_*$.  
Hence \eqref{eq-cor} and \eqref{eq-direction} hold for $c>c_*$.
\end{enumerate}
\end{exam}

\section{Proofs of Theorems \ref{thm-main} and \ref{thm-direction}}\label{sect-prf}

\subsection{Growth rate of Feynman-Kac semigroups}
\label{subsect-fk}

In this subsection, we reveal the growth rate of the expectation of $Z_t^{\delta t}$. 
We recall that 
$${\mathbf E}_x\left[Z_t^{\delta t}\right]=E_x\left[e^{A_t^{(Q-1)\mu}};|B_t|\geq \delta t\right]$$
by \eqref{eq-br-fk}.  
When $d=1$ and $\mu$ is the Dirac measure,
Bocharov and Harris \cite[Proposition 2]{BH14} calculated the right hand side above 
by using the joint distribution of the Brownian motion and local time. 
We first show a similar result when $\mu$ is a Kato class measure 
with compact support in ${\mathbb R}^d$ such that $\lambda(\mu)<0$.

\begin{prop}\label{prop-asymp}
Let $\mu$ be a Kato class measure with compact support in ${\mathbb R}^d$. 
If $\lambda(\mu)<0$, then for any $x\in {\mathbb R}^d$ and $\delta>0$, 
\begin{equation}\label{eq-fkf}
E_x\left[e^{A_t^{\mu}}; |B_t|\geq \delta t\right]\asymp 
\begin{cases} 
e^{\left(-\lambda(\mu)-\sqrt{-2\lambda(\mu)}\delta\right) t}t^{(d-1)/2} 
& \text{if $\delta\leq \sqrt{-2\lambda(\mu)}$,}\\
e^{-\delta^2 t/2}t^{(d-2)/2} & \text{if $\delta>\sqrt{-2\lambda(\mu)}$}
\end{cases}
\quad (t\rightarrow\infty).
\end{equation}
\end{prop}

Let us explain the meaning of Proposition \ref{prop-asymp}.
Since $\mu$ is compactly supported by assumption, we see that if $\delta$ is large enough, 
then the effect of the Feynman-Kac functional $e^{A_t^{\mu}}$ disappears 
at the left hand side of \eqref{eq-fkf}. In fact, 
\begin{equation*}
\begin{split}
E_x[e^{A_t^{\mu}};|B_t|\geq \delta t]
\asymp P_x\left(|B_t|\geq \delta t\right)
&\asymp \frac{1}{(2\pi t)^{d/2}}\int_{|y|\geq \delta t}e^{-|y|^2/(2t)}\,{\rm d}y \\
&\asymp e^{-\delta^2 t/2}t^{(d-2)/2} 
 \quad (t\rightarrow\infty).
\end{split}
\end{equation*}
The last relation above follows from \eqref{eq-int-comp} below. 

On the other hand, if $\delta$ is small enough, 
then the effect of the Feynman-Kac functional appears through  
the eigenfunction $h$. 
If $\delta=0$, then by taking $f(x)\equiv 1$ in \eqref{eq-fkf-asymp}, 
we obtain 
\begin{equation}\label{eq-feyn-lim}
\lim_{t\rightarrow\infty}
e^{\lambda(\mu) t}E_x\left[e^{A_t^{\mu}}\right]=c_h(x) 
\left(:=h(x)\int_{{\mathbb R}^d}h(y)\,{\rm d}y\right)
\quad \text{for any $x\in {\mathbb R}^d$}.
\end{equation}
Proposition \ref{prop-asymp} says that if $\delta\in (0,\sqrt{-2\lambda(\mu)}]$, then
\eqref{eq-fkf-asymp} is valid even for $f(x)={\bf 1}_{|x|\geq \delta t}$:
\begin{equation}\label{eq-asymp-time}
\begin{split}
E_x[e^{A_t^{\mu}};|B_t|\geq \delta t]
\asymp e^{-\lambda(\mu) t}\int_{|y|\geq \delta t}h(y)\,{\rm d}y
&\asymp e^{-\lambda(\mu) t}\int_{|y|\geq \delta t}\frac{e^{-\sqrt{-2\lambda}|y|}}{|y|^{(d-1)/2}}\,{\rm d}y\\
&\asymp e^{\left(-\lambda(\mu)-\sqrt{-2\lambda(\mu)}\delta\right) t}t^{(d-1)/2} 
\quad (t\rightarrow\infty).
\end{split}
\end{equation}
The second relation above follows from \eqref{eq-ground}. 
\medskip

\pf \ 
Let $\mu$ be a Kato class measure with compact support in ${\mathbb R}^d$. 
Define $A_t:=A_t^{\mu}$ and $\lambda:=\lambda(\mu)$. 
Since 
$$e^{A_t}=1+\int_0^t e^{A_s}\,{\rm d}A_s,$$
we have 
\begin{equation*}
\begin{split}
E_x\left[e^{A_t}; |B_t|\geq \delta t\right]
&=E_x\left[1+\int_0^t e^{A_s}\,{\rm d}A_s ; |B_t|\geq \delta t\right]\\
&=P_x(|B_t|\geq \delta t)+E_x\left[\int_0^t e^{A_s}{\bf 1}_{\{|B_t|\geq \delta t\}}\,{\rm d}A_s\right]
={\rm (I)}+{\rm (II)}.
\end{split}
\end{equation*}
Then by Appendix \ref{appendix-a-1} below,   
\begin{equation}\label{eq-lower-ball}
\begin{split}
{\rm (I)}
\geq P_0(|B_t|\geq \delta t)
=c\int_{\delta \sqrt{t}}^{\infty}\exp\left(-\frac{r^2}{2}\right)r^{d-1}\,{\rm d}r
\end{split}
\end{equation}
and 
\begin{equation}\label{eq-upper-ball}
\begin{split}
{\rm (I)}
=\frac{1}{(2\pi t)^{d/2}}\int_{|y|\geq \delta t}\exp\left(-\frac{|y-x|^2}{2t}\right)\,{\rm d}y
&\leq \frac{1}{(2\pi t)^{d/2}}\int_{|y-x|\geq \delta t-|x|}\exp\left(-\frac{|y-x|^2}{2t}\right)\,{\rm d}y\\
&=c\int_{(\delta t-|x|)/\sqrt{t}}^{\infty}\exp\left(-\frac{r^2}{2}\right)r^{d-1}\,{\rm d}r.
\end{split}
\end{equation}
Since
\begin{equation}\label{eq-int-comp}
\int_t^{\infty}\exp\left(-\frac{v^2}{2}\right)v^{d-1}\,{\rm d}v
\sim \exp\left(-\frac{t^2}{2}\right)t^{d-2} \quad (t\rightarrow\infty),
\end{equation}
we get
\begin{equation}\label{eq-i}
{\rm (I)}\sim c'\exp\left(-\frac{\delta^2 t}{2}\right)t^{(d-2)/2}  \quad (t\rightarrow\infty).
\end{equation}

We first discuss the lower bound of ${\rm (I)}+{\rm (II)}$. 
Fix $t>0$ and $R>0$. Let us define 
$$X_s=e^{A_s}P_x(|B_t|\geq R \mid {\cal F}_s){\bf 1}_{\{s\leq t\}}, \quad Y_s=e^{A_s}{\bf 1}_{\{|B_t|\geq R\}}{\bf 1}_{\{s\leq t\}} 
\quad (s\geq 0).$$
Since the strong Markov property implies that 
$$E_x[X_T{\bf 1}_{\{T<\infty\}}]=E_x\left[Y_T{\bf 1}_{\{T<\infty\}}\right]$$
for any $\{{\cal F}_t\}$-stopping time $T$, we have by \cite[p.186, Exercise 1.13]{RY99}, 
\begin{equation}\label{eq-markov}
{\rm (II)}
=E_x\left[\int_0^t e^{A_s}P_x(|B_t|\geq \delta t \mid {\cal F}_s)\,{\rm d}A_s\right]
=E_x\left[\int_0^t e^{A_s}P_{B_s}(|B_{t-s}|\geq \delta t)\,{\rm d}A_s\right].
\end{equation}
Noting that by Appendix \ref{appendix-a-1} below, 
\begin{equation}\label{eq-lower-ball-0}
P_x(|B_{t-s}|\geq \delta t)
\geq P_0(|B_{t-s}|\geq \delta t)
=c \int_{\delta t/\sqrt{t-s}}^{\infty}\exp\left(-\frac{r^2}{2}\right)r^{d-1}\,{\rm d}r,
\end{equation}
we obtain 
$$({\rm II})
\geq c
E_x\left[\int_0^t e^{A_s}\left(\int_{\delta t/\sqrt{t-s}}^{\infty}\exp\left(-\frac{r^2}{2}\right)r^{d-1}\,{\rm d}r\right)\,{\rm d}A_s\right].$$
Then by the integration by parts formula, 
\begin{equation*}
\begin{split}
&\int_0^t e^{A_s}\left(\int_{\delta t/\sqrt{t-s}}^{\infty}\exp\left(-\frac{r^2}{2}\right)r^{d-1}\,{\rm d}r\right)\,{\rm d}A_s\\
&=\left[e^{A_s}\left(\int_{\delta t/\sqrt{t-s}}^{\infty}\exp\left(-\frac{r^2}{2}\right)r^{d-1}\,{\rm d}r\right)\right]_{s=0}^{s=t}
-\int_0^t e^{A_s}
\left[\frac{\partial}{\partial s}
\left(\int_{\delta t/\sqrt{t-s}}^{\infty}\exp\left(-\frac{r^2}{2}\right)r^{d-1}\,{\rm d}r\right)\right]\,{\rm d}s\\
&=-\int_{\delta \sqrt{t}}^{\infty}\exp\left(-\frac{r^2}{2}\right)r^{d-1}\,{\rm d}r
+\frac{(\delta t)^d}{2}\int_0^t e^{A_s}\exp \left(-\frac{(\delta t)^2}{2(t-s)}\right)\frac{1}{(t-s)^{(d+2)/2}}\,{\rm d}s.
\end{split}
\end{equation*}
Combining this with \eqref{eq-lower-ball}, we have 
$$
({\rm II})
\geq -({\rm I})
+c\frac{(\delta t)^d}{2}\int_0^t E_x\left[e^{A_s}\right]\exp\left(-\frac{(\delta t)^2}{2(t-s)}\right)
\frac{1}{(t-s)^{(d+2)/2}}\,{\rm d}s,
$$
that is, 
\begin{equation}\label{eq-iii-0}
({\rm I})+({\rm II})
\geq 
c\frac{(\delta t)^d}{2}\int_0^t E_x\left[e^{A_s}\right]\exp\left(-\frac{(\delta t)^2}{2(t-s)}\right)
\frac{1}{(t-s)^{(d+2)/2}}\,{\rm d}s
=\frac{c \delta^d}{2}({\rm III}).
\end{equation}

Since $\lambda<0$ by assumption, 
\eqref{eq-feyn-lim} yields that 
for any $x\in {\mathbb R}^d$ and $\varepsilon>0$, there exists $T>0$ such that 
\begin{equation}\label{eq-comp-funct}
(c_h(x)-\varepsilon)e^{-\lambda t}\leq 
E_x[e^{A_t}]\leq (c_h(x)+\varepsilon)e^{-\lambda t}
\quad \text{for all $t\geq T$}.
\end{equation}
For $t\geq T$, we let 
\begin{equation}\label{eq-iii-split}
\begin{split}
({\rm III})
&=t^d\int_0^T E_x\left[e^{A_s}\right]\exp \left(-\frac{(\delta t)^2}{2(t-s)}\right)\frac{1}{(t-s)^{(d+2)/2}}\,{\rm d}s\\
&+t^d\int_T^t E_x\left[e^{A_s}\right]\exp \left(-\frac{(\delta t)^2}{2(t-s)}\right)\frac{1}{(t-s)^{(d+2)/2}}\,{\rm d}s\\
&=({\rm III})_1+({\rm III})_2.
\end{split}
\end{equation}
Since 
\begin{equation*}
\begin{split}
\exp \left(-\frac{(\delta t)^2}{2(t-T)}\right)\frac{1}{t^{(d+2)/2}}
&\leq \exp \left(-\frac{(\delta t)^2}{2(t-s)}\right)\frac{1}{(t-s)^{(d+2)/2}}\\
&\leq \exp \left(-\frac{\delta^2t}{2}\right)\frac{1}{(t-T)^{(d+2)/2}}
\end{split}
\end{equation*}
for any $s\in [0,T]$ and $\sup_{x\in {\mathbb R}^d}E_x[e^{A_T}]<\infty$ (see, e.g., \cite[Theorem 6.1]{ABM91}), 
we obtain 
\begin{equation}\label{eq-a}
({\rm III})_1
\asymp t^{(d-2)/2}\exp\left(-\frac{\delta^2 t}{2}\right) 
\int_0^T E_x[e^{A_s}]\,{\rm d}s
\quad (t\rightarrow\infty).
\end{equation} 
By \eqref{eq-comp-funct}, we also get 
\begin{equation}\label{eq-b}
({\rm III})_2\asymp t^d\int_T^t e^{-\lambda s}\exp \left(-\frac{(\delta t)^2}{2(t-s)}\right)\frac{1}{(t-s)^{(d+2)/2}}\,{\rm d}s 
\quad (t\rightarrow\infty).
\end{equation}
By the change of variables $u=t-s$, 
\begin{equation}\label{eq-large}
\begin{split}
&\int_T^t e^{-\lambda s}\exp \left(-\frac{(\delta t)^2}{2(t-s)}\right)\frac{1}{(t-s)^{(d+2)/2}}\,{\rm d}s
=\int_0^{t-T} e^{-\lambda (t-u)}\exp\left(-\frac{(\delta t)^2}{2u}\right)\frac{1}{u^{(d+2)/2}}\,{\rm d}u\\
&=e^{(-\lambda-\sqrt{-2\lambda}\delta) t}\int_0^{t-T}
e^{-(\sqrt{-\lambda u}-\delta t/\sqrt{2u})^2}\frac{1}{u^{(d+2)/2}}\,{\rm d}u.
\end{split}
\end{equation}
If we further let $v=\sqrt{-\lambda u}-\delta t/\sqrt{2u}$, that is,
$$\sqrt{u}=\frac{v+\sqrt{v^2+2\sqrt{-2\lambda}\delta t}}{2\sqrt{-\lambda}},$$
then
$$\frac{{\rm d}u}{\sqrt{u}}
=\frac{v+\sqrt{v^2+2\sqrt{-2\lambda}\delta t}}{\sqrt{-\lambda}\sqrt{v^2+2\sqrt{-2\lambda}\delta t}}\,{\rm d}v.
$$
Hence \eqref{eq-large} is equal to 
$$2e^{(-\lambda-\sqrt{-2\lambda}\delta) t} \int_{-\infty}^S e^{-v^2}
\left(\frac{2\sqrt{-\lambda}}{v+\sqrt{v^2+2\sqrt{-2\lambda}\delta t}}\right)^d
\frac{1}{\sqrt{v^2+2\sqrt{-2\lambda}\delta t}}\,{\rm d}v$$
for
$$S=\sqrt{-\lambda(t-T)}-\frac{\delta t}{\sqrt{2(t-T)}}.$$
Namely, we have by \eqref{eq-b},
\begin{equation}\label{eq-c}
\begin{split}
({\rm III})_2
&\asymp 
t^d e^{(-\lambda-\sqrt{-2\lambda}\delta) t} 
\int_{-\infty}^S e^{-v^2}
\left(\frac{2\sqrt{-\lambda}}{v+\sqrt{v^2+2\sqrt{-2\lambda}\delta t}}\right)^d
\frac{1}{\sqrt{v^2+2\sqrt{-2\lambda}\delta t}}\,{\rm d}v\\
&=t^d e^{(-\lambda-\sqrt{-2\lambda}\delta) t} ({\rm IV)} \quad (t\rightarrow\infty).
\end{split}
\end{equation}

We will prove in Appendix \ref{appendix-a-2} below that 
\begin{equation}\label{eq-iv}
({\rm IV})\asymp 
\begin{cases}
t^{-(d+1)/2} & \text{if $\delta\leq \sqrt{-2\lambda}$},\\
 e^{-(-\lambda-\sqrt{-2\lambda}\delta) t}e^{-\delta^2t/2}t^{-(d+2)/2} & 
\text{if $\delta>\sqrt{-2\lambda}$}
\end{cases}
\quad (t\rightarrow\infty).
\end{equation}
Since this relation implies that 
$$
({\rm III})_2
\asymp
\begin{cases}
e^{(-\lambda-\sqrt{-2\lambda}\delta)t}t^{(d-1)/2}  & \text{if $\delta\leq \sqrt{-2\lambda}$},\\
e^{-\delta^2t/2}t^{(d-2)/2} & \text{if $\delta> \sqrt{-2\lambda}$}
\end{cases}
\quad (t\rightarrow\infty),$$
we obtain by \eqref{eq-iii-0}, \eqref{eq-iii-split} and \eqref{eq-a},
\begin{equation}\label{eq-iii}
({\rm I})+({\rm II})\geq \frac{c\delta^d}{2}\left(({\rm III})_1+({\rm III})_2\right)
\asymp
\begin{cases}
e^{(-\lambda-\sqrt{-2\lambda}\delta)t}t^{(d-1)/2}  & \text{if $\delta\leq \sqrt{-2\lambda}$},\\
e^{-\delta^2t/2}t^{(d-2)/2} & \text{if $\delta>\sqrt{-2\lambda}$}
\end{cases}
\quad (t\rightarrow\infty).
\end{equation}

We next discuss the upper bound of $({\rm I})+({\rm II})$.
Here we suppose that the support of $\mu$ is included 
in a ball $\{x\in {\mathbb R}^d \mid |x|\leq R\}$ for some $R>0$.
Then by the same calculation as in \eqref{eq-upper-ball}, 
we see from \eqref{eq-markov} that 
\begin{equation}\label{eq-ii-upper}
\begin{split}
({\rm II})
&\leq cE_x\left[\int_0^t e^{A_s}
\left(
\int_{(\delta t-|B_s|)/\sqrt{t-s}}^{\infty}\exp\left(-\frac{r^2}{2}\right)r^{d-1}\,{\rm d}r\right)
\,{\rm d}A_s\right]\\
&\leq c
E_x\left[\int_0^t e^{A_s}
\left(\int_{(\delta t-R)/\sqrt{t-s}}^{\infty}\exp\left(-\frac{r^2}{2}\right)r^{d-1}\,{\rm d}r\right)
\,{\rm d}A_s\right].
\end{split}
\end{equation}
By the integration by parts formula, 
\begin{equation*}
\begin{split}
&\int_0^t e^{A_s}\left(\int_{(\delta t-R)/\sqrt{t-s}}^{\infty}
\exp\left(-\frac{r^2}{2}\right)r^{d-1}\,{\rm d}r\right)\,{\rm d}A_s\\
&=-\int_{(\delta t-R)/\sqrt{t}}^{\infty}
\exp\left(-\frac{r^2}{2}\right)r^{d-1}\,{\rm d}r
+\frac{(\delta t-R)^d}{2}\int_0^t e^{A_s}\exp \left(-\frac{(\delta t-R)^2}{2(t-s)}\right)\frac{1}{(t-s)^{(d+2)/2}}\,{\rm d}s\\
&\leq 
\frac{(\delta t-R)^d}{2}\int_0^t e^{A_s}\exp \left(-\frac{(\delta t-R)^2}{2(t-s)}\right)\frac{1}{(t-s)^{(d+2)/2}}\,{\rm d}s,
\end{split}
\end{equation*}
whence
$$
({\rm II})
\leq c
\frac{(\delta t-R)^d}{2}\int_0^t E_x[e^{A_s}]
\exp \left(-\frac{(\delta t-R)^2}{2(t-s)}\right)\frac{1}{(t-s)^{(d+2)/2}}\,{\rm d}s.
$$
Then  
$$
\frac{(\delta t-R)^d}{2}\int_0^T E_x\left[e^{A_s}\right]
\exp \left(-\frac{(\delta t-R)^2}{2(t-s)}\right)\frac{1}{(t-s)^{(d+2)/2}}\,{\rm d}s
\asymp e^{-\delta^2 t/2}t^{(d-2)/2} \quad (t\rightarrow\infty)
$$
by the same way as for the evaluation of  $({\rm III})_1$.
We also have by \eqref{eq-comp-funct},
\begin{equation*}
\begin{split}
&\frac{(\delta t-R)^d}{2}\int_T^t E_x\left[e^{A_s}\right]
\exp \left(-\frac{(\delta t-R)^2}{2(t-s)}\right)\frac{1}{(t-s)^{(d+2)/2}}\,{\rm d}s\\
&\asymp  
\frac{(\delta t-R)^d}{2}\int_T^t e^{-\lambda s}\exp \left(-\frac{(\delta t-R)^2}{2(t-s)}\right)\frac{1}{(t-s)^{(d+2)/2}}\,{\rm d}s 
\quad (t\rightarrow\infty).
\end{split}
\end{equation*}
If we let 
$$S'=\sqrt{-\lambda(t-T)}-\frac{\delta t-R}{\sqrt{2(t-T)}},$$
then 
$$S'=\left(\sqrt{-\lambda}-\frac{\delta}{\sqrt{2}}-\frac{\delta T-R}{\sqrt{2}(t-T)}\right)\sqrt{t-T}.$$
Hence if we take  $T>0$ so large that $\delta T>R$, 
then the next relation follows  
by the same argument as for the evaluation of $({\rm III})_2$:
\begin{equation*}
\begin{split}
&\frac{(\delta t-R)^d}{2}\int_T^t e^{-\lambda s}
\exp \left(-\frac{(\delta t-R)^2}{2(t-s)}\right)\frac{1}{(t-s)^{(d+2)/2}}\,{\rm d}s \\
&\asymp
\begin{cases} 
e^{(-\lambda-\sqrt{-2\lambda}\delta)t}t^{(d-1)/2} & 
\text{if $\delta\leq\sqrt{-2\lambda}$},\\
e^{-\delta^2t/2}t^{(d-2)/2} & \text{if $\delta>\sqrt{-2\lambda}$}
\end{cases} 
\quad (t\rightarrow\infty).
\end{split}
\end{equation*}
Combining this with \eqref{eq-i}, we get
\begin{equation}\label{eq-upper-i-ii}
({\rm I})+({\rm II})\lesssim
\begin{cases}
e^{(-\lambda-\sqrt{-2\lambda}\delta)t}t^{(d-1)/2}  & \text{if $\delta\leq \sqrt{-2\lambda}$},\\
e^{-\delta^2t/2}t^{(d-2)/2} & \text{if $\delta>\sqrt{-2\lambda}$}
\end{cases}
\quad (t\rightarrow\infty),
\end{equation}
which completes the proof.
\qed

\begin{rem}\label{rem-asymp} \rm
As we see from the proof of Proposition \ref{prop-asymp}, 
the lower bound of ${\rm (I)}+{\rm (II)}$ is true 
for any $\mu\in {\cal K}_{\infty}(1)$. 
\end{rem}

Our next result says that 
we can relax the compact support condition on $\mu$ in Proposition \ref{prop-asymp} 
only to find the exponential growth rate of 
$E_x\left[e^{A_t^{\mu}}; |B_t|\geq \delta t\right]$.

\begin{prop}\label{prop-non-cpt}
Let $\mu$ be a positive Radon measure on ${\mathbb R}^d$ such that 
the measure $\nu_{\beta}({\rm d}x):=e^{\beta |x|}\mu({\rm d}x)$ 
belongs to ${\cal K}_{\infty}(1)$ for any $\beta>0$. 
If $\lambda(\mu)<0$, then for any $\delta\geq 0$,
$$\lim_{t\rightarrow\infty}\frac{1}{t}\log E_x\left[e^{A_t^{\mu}}; |B_t|\geq \delta t\right]
=\begin{cases} 
-\lambda(\mu)-\sqrt{-2\lambda(\mu)}\delta& 
\text{if $\delta\leq \sqrt{-2\lambda(\mu)}$},\\
-\delta^2/2 & \text{if $\delta>\sqrt{-2\lambda(\mu)}$.}
\end{cases}$$
\end{prop}

The condition on $\nu_{\beta}$ means that the measure $\mu$ is small enough at infinity. 
This condition is a key to show that 
the first expression in \eqref{eq-non-cpt-add} below is uniformly bounded.  
\smallskip

\pf \ 
If $\delta=0$, then our assertion follows from \eqref{eq-feyn-lim}. 
In what follows, we assume that $\delta>0$. 
Let $\mu$ be a positive Radon measure on ${\mathbb R}^d$ 
satisfying the condition in this proposition.
Define $A_t:=A_t^{\mu}$ and $\lambda:=\lambda(\mu)$. 
Since $\mu\in {\cal K}_{\infty}(1)$,  
we see by Remark \ref{rem-asymp} (ii) that for any $\delta>0$,
$$\liminf_{t\rightarrow\infty}\frac{1}{t}\log E_x\left[e^{A_t}; |B_t|\geq \delta t\right]
\geq 
\begin{cases} 
-\lambda-\sqrt{-2\lambda}\delta& 
\text{if $\delta\leq \sqrt{-2\lambda}$},\\
-\delta^2/2 & \text{if $\delta>\sqrt{-2\lambda}$.}
\end{cases}
$$

We next show that 
\begin{equation}\label{eq-limsup}
\limsup_{t\rightarrow\infty}
\frac{1}{t}\log E_x\left[e^{A_t}; |B_t|\geq \delta t\right]
\leq 
\begin{cases} 
-\lambda-\sqrt{-2\lambda}\delta& 
\text{if $\delta\leq \sqrt{-2\lambda}$},\\
-\delta^2/2 & \text{if $\delta>\sqrt{-2\lambda}$}
\end{cases}
\end{equation}
by modifying the argument 
for the upper bound of ${\rm (I)}+{\rm (II)}$ in the proof of Proposition \ref{prop-asymp}. 
We keep the same notations as in the proof of Proposition \ref{prop-asymp}. 
For any $\varepsilon\in (0,\delta)$, we have by \eqref{eq-ii-upper}, 
\begin{equation*}
\begin{split}
({\rm II})
&\leq cE_x\left[\int_0^t e^{A_s}
\left(
\int_{(\delta t-|B_s|)/\sqrt{t-s}}^{\infty}\exp\left(-\frac{r^2}{2}\right)r^{d-1}\,{\rm d}r\right)
\,{\rm d}A_s \right]\\
&\leq c
E_x\left[\int_0^t e^{A_s}
\left(\int_{(\delta -{\varepsilon}) t/\sqrt{t-s}}^{\infty}\exp\left(-\frac{r^2}{2}\right)r^{d-1}\,{\rm d}r\right)
{\bf 1}_{\{|B_s|< \varepsilon t\}}\,{\rm d}A_s\right]\\
&+
E_x\left[\int_0^t e^{A_s}
{\bf 1}_{\{|B_s|\geq \varepsilon t\}}\,{\rm d}A_s\right].
\end{split}
\end{equation*}
Then by the same way as in the proof of Proposition \ref{prop-asymp}, 
\begin{equation*}
\begin{split}
&E_x\left[\int_0^t e^{A_s}
\left(\int_{(\delta -{\varepsilon}) t/\sqrt{t-s}}^{\infty}\exp\left(-\frac{r^2}{2}\right)r^{d-1}\,{\rm d}r\right)
{\bf 1}_{\{|B_s|< \varepsilon t\}}\,{\rm d}A_s\right] \\
&\leq 
E_x\left[\int_0^t e^{A_s}
\left(\int_{(\delta -{\varepsilon}) t/\sqrt{t-s}}^{\infty}\exp\left(-\frac{r^2}{2}\right)r^{d-1}\,{\rm d}r\right)
\,{\rm d}A_s\right]\\
&\asymp
\begin{cases} 
e^{(-\lambda-\sqrt{-2\lambda}(\delta-\varepsilon))t}t^{(d-1)/2} & 
\text{if $\delta\leq \sqrt{-2\lambda}+\varepsilon$},\\
e^{-(\delta-\varepsilon)^2t/2}t^{(d-2)/2} 
& \text{if $\delta>\sqrt{-2\lambda}+\varepsilon$}
\end{cases}
\quad (t\rightarrow\infty).
\end{split}
\end{equation*}

For any $c>0$,
\begin{equation}\label{eq-non-cpt-add}
\begin{split}
E_x\left[\int_0^t e^{A_s}
{\bf 1}_{\{|B_s|\geq \varepsilon t\}}\,{\rm d}A_s\right]
&\leq 
e^{-c t}E_x\left[\int_0^{\infty} e^{-cs}e^{A_s}
{\bf 1}_{\{|B_s|\geq \varepsilon t\}}e^{2c\varepsilon^{-1}|B_s|}\,{\rm d}A_s\right]\\
&\leq e^{-c t}E_x\left[\int_0^{\infty} e^{-cs}e^{A_s}
\,{\rm d}A_s^{\nu_{2c/\varepsilon}}\right].
\end{split}
\end{equation}
If we take $c>0$ so large that $c>-\lambda$, then 
$$\inf\left\{\frac{1}{2}\int_{{\mathbb R}^d}|\nabla u|^2\,{\rm d}x
+c\int_{{\mathbb R}^d}u^2\,{\rm d}x-\int_{{\mathbb R}^d}u^2\,{\rm d}\mu
\mid u\in C_0^{\infty}({\mathbb R}^d), \int_{{\mathbb R}^d}u^2\,{\rm d}x=1\right\}
=c+\lambda>0$$
and thus 
$$
\inf\left\{\frac{1}{2}\int_{{\mathbb R}^d}|\nabla u|^2\,{\rm d}x
+c\int_{{\mathbb R}^d}u^2\,{\rm d}x
\mid u\in C_0^{\infty}({\mathbb R}^d), \int_{{\mathbb R}^d}u^2\,{\rm d}\mu=1\right\}
>1$$
by \cite[Lemma 3.5]{T02}. 
Since $\nu_{2c/\varepsilon}\in {\cal K}_{\infty}(1)$ by assumption, 
we get 
$$
\sup_{x\in {\mathbb R}^d}
E_x\left[\int_0^{\infty} e^{-c s}e^{A_s}
\,{\rm d}A_s^{\nu_{2c/\varepsilon}}\right]<\infty$$
by \cite[Corollary 2.9 and Theorem 5.1]{C02}.
Hence by taking $c>0$ large enough, we see that 
$${\rm (II)}
\lesssim 
\begin{cases} 
e^{(-\lambda-\sqrt{-2\lambda}(\delta-\varepsilon))t}t^{(d-1)/2} & 
\text{if $\delta\leq \sqrt{-2\lambda}+\varepsilon$},\\
e^{-(\delta-\varepsilon)^2t/2}t^{(d-2)/2} 
& \text{if $\delta>\sqrt{-2\lambda}+\varepsilon$}
\end{cases} 
\quad (t\rightarrow\infty),$$
that is,
$$\limsup_{t\rightarrow\infty}\frac{1}{t}\log({\rm (I)}+{\rm (II)})
\leq 
\begin{cases} 
-\lambda-\sqrt{-2\lambda}(\delta-\varepsilon)& 
\text{if $\delta\leq\sqrt{-2\lambda}+\varepsilon$},\\
-(\delta-\varepsilon)^2/2& \text{if $\delta>\sqrt{-2\lambda}+\varepsilon$.}
\end{cases}$$
By letting $\varepsilon\rightarrow+0$, we arrive at  \eqref{eq-limsup}. 
\qed

\begin{rem}\rm 
Let $r$ be a unit vector in ${\mathbb R}^d$ 
and $\langle \cdot,\cdot\rangle$ the standard inner product on ${\mathbb R}^d$.
Then Proposition  \ref{prop-asymp} is
still true by replacing $|B_t|$ with $\langle B_t, r\rangle$, 
and 
by taking $d$ as $1$ at the right hand side of \eqref{eq-fkf}. 
For $d=1$, we can replace $|B_t|$ with $B_t$.   
These modifications are also applicable to Proposition \ref{prop-non-cpt}. 
We can prove these assertions in a similar way to Propositions \ref{prop-asymp} and \ref{prop-non-cpt}. 

Recall that 
$B_R^r=\{x\in {\mathbb R}^d \mid \langle x, r\rangle <R\}$ for $R>0$ and 
$Z_t^{\delta t,r}=Z_t((B_{\delta t}^r)^c)$.
Since 
$${\mathbf E}_x\left[Z_t^{\delta t,r}\right]=E_x\left[e^{A_t^{(Q-1)\mu}};\langle B_t,r\rangle\geq \delta t\right]$$
by \eqref{eq-br-fk}, 
the modifications of Propositions \ref{prop-asymp} and \ref{prop-non-cpt} 
show that the growth rate of ${\mathbf E}_x\left[Z_t^{\delta t,r}\right]$ is uniform in direction. 
\end{rem}

Let $Q_x^{\mu,(t)}$ be the probability measure on $(\Omega, {\cal F}_t)$ defined by 
$$Q_x^{\mu, (t)}({\rm d}\omega):=\frac{e^{A_t^{\mu}(\omega)}P_x({\rm d}\omega)}{E_x[e^{A_t^{\mu}}]}.$$
This measure gives the distribution of the Brownian motion 
under the probability measure weighted by the Feynman-Kac functional. 
Then by Proposition \ref{prop-non-cpt}, 
we obtain the large deviation type result as follows:

\begin{cor}
Under the same setting as in Proposition {\rm \ref{prop-non-cpt}}, 
if $\lambda(\mu)<0$, then for any $\delta>0$,
$$\lim_{t\rightarrow\infty}\frac{1}{t}\log Q_x^{\mu,(t)}(|B_t|\geq \delta t)
=\begin{cases} 
-\sqrt{-2\lambda(\mu)}\delta& 
\text{if $\delta\leq \sqrt{-2\lambda(\mu)}$},\\
\lambda(\mu)-\delta^2/2 & \text{if $\delta>\sqrt{-2\lambda(\mu)}$.}
\end{cases}$$
\end{cor}

\subsection{Preliminary lemmas}

We now discuss the expected number of particles 
conditioned on their trajectories. 
For a fixed $t>0$, we denote by ${\mathbf B}_s^{(t),k}$ 
the position of the ancestor of ${\mathbf B}_t^k$ at time $s$ 
($s\leq t$ and $1\leq k\leq Z_t$).

\begin{lem}\label{lem-exp}
Suppose that $\nu$ is a Kato class measure. 
Then for any $\delta>0$ and $t>0$, 
$${\mathbf E}_x\left[
\sum_{k=1}^{Z_t}{\bf 1}_{\left\{\sup_{0\leq s\leq t}|{\mathbf B}_s^{(t),k}|\geq \delta t\right\}}\right]
=E_x\left[e^{A_t^{(Q-1)\mu}} ; \sup_{0\leq s\leq t}|B_s|\geq \delta t\right].$$
\end{lem}

\pf \ We prove this lemma in a similar way to \eqref{eq-br-fk} 
(see \cite[Lemma 3.3]{S08}).
Let $Z_t(n)$ be the total number of particles at time $t$ such that 
each of their trajectories over the time interval $[0,t]$ has $n$ branching points, 
and 
$${\mathbf B_t}(n)=({\mathbf B}_t^1(n),\dots, {\mathbf B}_t^{Z_t(n)}(n))$$ 
the positions of such particles at time $t$. 
Denote by ${\mathbf B}_s^{(t),k}(n)$ 
the position of the ancestor of ${\mathbf B}_t^k(n)$ at time $s$. 
We will prove that 
\begin{equation}\label{eq-n-pt}
{\mathbf E}_x\left[
\sum_{k=1}^{Z_t(n)}
{\bf 1}_{\left\{\sup_{0\leq s\leq t}|{\mathbf B}_s^{(t),k}(n)|\geq \delta t\right\}}\right]
=E_x\left[e^{-A_t^{\mu}}\frac{(A_t^{Q\mu})^n}{n!}; \sup_{0\leq s\leq t}|B_s|\geq \delta t\right] 
\quad \text{for any $n\geq 0$}.
\end{equation}
Using this equality, we obtain 
\begin{equation*}
\begin{split}
&{\mathbf E}_x\left[
\sum_{k=1}^{Z_t}
{\bf 1}_{\left\{\sup_{0\leq s\leq t}|{\mathbf B}_s^{(t),k}|\geq \delta t\right\}}\right]
=\sum_{n=0}^{\infty}{\mathbf E}_x\left[
\sum_{k=1}^{Z_t(n)}
{\bf 1}_{\left\{\sup_{0\leq s\leq t}|{\mathbf B}_s^{(t),k}(n)|\geq \delta t\right\}}\right]\\
&=\sum_{n=0}^{\infty}E_x\left[e^{-A_t^{\mu}}\frac{(A_t^{Q\mu})^n}{n!} ; \sup_{0\leq s\leq t}|B_s|\geq \delta t\right]
=E_x\left[e^{A_t^{(Q-1)\mu}}; \sup_{0\leq s\leq t}|B_s|\geq \delta t\right],
\end{split}
\end{equation*}
which is our assertion.

In what follows, we show \eqref{eq-n-pt} by induction. 
For $n=0$, we have by \eqref{eq-split},
\begin{equation*}\label{eq-non-branching}
\begin{split}
{\mathbf E}_x\left[
\sum_{k=1}^{Z_t(0)}
{\bf 1}_{\left\{\sup_{0\leq s\leq t}|{\mathbf B}_s^{(t),k}(0)|\geq \delta t\right\}}\right]
&={\mathbf P}_x
\left(t<T,  \sup_{0\leq s\leq t} |{\mathbf B}_s^1|\geq \delta t \right)\\
&=E_x\left[e^{-A_t^{\mu}}; \sup_{0\leq s\leq t} |B_s|\geq \delta t\right].
\end{split}
\end{equation*}
Suppose that \eqref{eq-n-pt} is valid for some $n\geq 1$. 
By the strong Markov property, 
\begin{equation}\label{eq-st}
\begin{split}
&{\mathbf E}_x
\left[
\sum_{k=1}^{Z_t(n+1)}
{\bf 1}_{\left\{\sup_{0\leq s\leq t}|{\mathbf B}_s^{(t),k}(n+1)|\geq \delta t\right\}}\right]\\
&={\mathbf E}_x
\left[{\mathbf E}_{{\mathbf B}_T}
\left[Z_{t-s}(n)\right]|_{s=T}
; \sup_{0\leq u\leq T}|{\mathbf B}_u^1|\geq \delta t, T\leq t\right]\\
&+{\mathbf E}_x\left[{\mathbf E}_{{\mathbf B}_T}
\left[\sum_{k=1}^{Z_{t-s}(n)}
{\bf 1}_{\left\{\sup_{0\leq u\leq t-s}|{\mathbf B}_u^{(t-s),k}(n)|\geq \delta t\right\}}
\right]|_{s=T}; 
\sup_{0\leq u\leq T}|{\mathbf B}_u^1|< \delta t, T\leq t\right]\\
&={\mathbf E}_x
\left[\sum_{k=1}^{Z_T}{\mathbf E}_{{\mathbf B}_T^k}
\left[Z_{t-s}(n)\right]|_{s=T}
; \sup_{0\leq u\leq T}|{\mathbf B}_u^1|\geq \delta t, T\leq t\right]\\
&+{\mathbf E}_x\left[\sum_{l=1}^{Z_T}{\mathbf E}_{{\mathbf B}_T^l}
\left[
\sum_{k=1}^{Z_{t-s}(n)}
{\bf 1}_{\left\{\sup_{0\leq u\leq t-s}|{\mathbf B}_u^{(t-s),k}(n)|\geq \delta t\right\}}
\right]|_{s=T}; 
\sup_{0\leq u\leq T}|{\mathbf B}_u^1|< \delta t, T\leq t\right].
\end{split}
\end{equation}
Then by \eqref{eq-n-pt} and 
$${\mathbf E}_x\left[Z_t(n)\right]
=E_x\left[e^{-A_t^{\mu}}\frac{(A_t^{Q\mu})^n}{n!}\right] \quad \text{for any $n\geq 0$}$$
(see, e.g., \cite[Lemma 3.3]{S08}), 
the last expression of \eqref{eq-st} is equal to 
\begin{equation*}
\begin{split}
&E_x\left[\int_0^t e^{-A_s^{\mu}}E_{B_s}\left[e^{-A_{t-s}^{\mu}}\frac{(A_{t-s}^{Q\mu})^n}{n!}\right]
{\bf 1}_{\{\sup_{0\leq u\leq s}|B_u|\geq \delta t\}}\,{\rm d}A_s^{Q\mu}\right]\\
&+E_x\left[\int_0^t e^{-A_s^{\mu}}E_{B_s}\left[e^{-A_{t-s}^{\mu}}\frac{(A_{t-s}^{Q\mu})^n}{n!}; 
\sup_{0\leq u\leq t-s}|B_u|\geq \delta t\right]
{\bf 1}_{\{\sup_{0\leq u\leq s}|B_u|< \delta t\}}\,{\rm d}A_s^{Q\mu}\right]\\
&=E_x\left[\int_0^t e^{-A_s^{\mu}}
E_x\left[e^{-A_{t-s}^{\mu}\circ\theta_s}\frac{(A_{t-s}^{Q\mu}\circ\theta_s)^n}{n!}; 
\sup_{0\leq u\leq t}|B_u|\geq \delta t \mid {\cal F}_s\right]\,{\rm d}A_s^{Q\mu}\right]\\
&=E_x\left[\int_0^t e^{-A_s^{\mu}}
e^{-A_{t-s}^{\mu}\circ\theta_s}\frac{(A_{t-s}^{Q\mu}\circ\theta_s)^n}{n!}
{\bf 1}_{\{\sup_{0\leq u\leq t}|B_u|\geq \delta t\}}\,{\rm d}A_s^{Q\mu}\right].
\end{split}
\end{equation*}
The last equality above follows by the same way as in the proof of \eqref{eq-markov}. 
By noting that $A_t^{\mu}$ and $A_t^{Q\mu}$ are additive functionals, 
the last expression above is equal to 
\begin{equation*}
\begin{split}
&E_x\left[e^{-A_t^{\mu}}\int_0^t\frac{(A_t^{Q\mu}-A_s^{Q\mu})^n}{n!}\,{\rm d}A_s^{Q\mu}
; \sup_{0\leq u\leq t}|B_u|\geq \delta t\right]\\
&=E_x\left[e^{-A_t^{\mu}}\frac{(A_t^{Q\mu})^{n+1}}{(n+1)!}; \sup_{0\leq u\leq t}|B_u|\geq \delta t\right],
\end{split}
\end{equation*}
which yields \eqref{eq-n-pt}. 
\qed 
\medskip

We next discuss the exponential growth rate of $Z_t$. 
\begin{lem}\label{lem-growth}
Suppose that $\nu\in {\cal K}_{\infty}(1)$. 
If $\lambda<0$ and ${\mathbf P}_x(M_{\infty}>0)>0$, then 
$$\lim_{t\rightarrow\infty}\frac{1}{t}\log Z_t
=\lim_{t\rightarrow\infty}\frac{1}{t}\log {\mathbf E}_x[Z_t]
=-\lambda \quad \text{${\mathbf P}_x(\cdot \mid M_{\infty}>0)$-a.s.}$$
\end{lem}

\pf \ The second equality above follows from \eqref{eq-fkf-asymp} and \eqref{eq-br-fk}. 
Since the eigenfunction $h$ is bounded (e.g., see \cite{T08}), 
we have 
$$M_t\leq e^{\lambda t}Z_t\|h\|_{\infty}$$
so that 
$$\liminf_{t\rightarrow\infty}\frac{1}{t}\log Z_t\geq -\lambda \quad 
\text{${\mathbf P}_x(\cdot \mid M_{\infty}>0)$-a.s.}$$

By Chebyshev's inequality, \eqref{eq-br-fk} and  \eqref{eq-comp-funct}, 
we obtain for any $\varepsilon>0$, 
$$
{\mathbf P}_x(e^{(\lambda-\varepsilon)t}Z_t>\varepsilon)
\leq \frac{e^{(\lambda-\varepsilon)t}}{\varepsilon}{\mathbf E}_x[Z_t]
=\frac{e^{(\lambda-\varepsilon)t}}{\varepsilon}E_x\left[e^{A_t^{(Q-1)\mu}}\right]
\leq \frac{c}{\varepsilon}e^{-\varepsilon t}
$$
for all sufficiently large $t>0$, whence
$$\sum_{n=1}^{\infty}{\mathbf P}_x(e^{(\lambda-\varepsilon)n}Z_n>\varepsilon)
\leq \frac{c}{\varepsilon}\sum_{n=1}^{\infty}e^{-\varepsilon n}<\infty.$$
Then by Borel-Cantelli's lemma, 
$${\mathbf P}_x\left(\text{$e^{(\lambda-\varepsilon)n}Z_n\leq \varepsilon$ 
for all sufficiently large $n$}\right)=1.$$
Since $Z_t$ is nondecreasing by assumption, we get 
$$Z_t\leq Z_{[t]+1}\leq \varepsilon e^{(-\lambda+\varepsilon)([t]+1)}
\leq \varepsilon e^{(-\lambda+\varepsilon)(t+1)}$$
for all sufficiently large $t>0$. Therefore,
$$\limsup_{t\rightarrow\infty}\frac{1}{t}\log Z_t\leq -\lambda+\varepsilon \quad \text{${\mathbf P}_x$-a.s.}$$
By letting $\varepsilon\rightarrow +0$, we obtain  
$$\limsup_{t\rightarrow\infty}\frac{1}{t}\log Z_t\leq -\lambda \quad \text{${\mathbf P}_x$-a.s.}$$
so that  the proof is complete. 
\qed
\medskip

\subsection{Proof of Theorem \ref{thm-main}}

Define for $\delta>0$, 
$$
\Lambda_{\delta}:=
\begin{cases}
\lambda+\sqrt{-2\lambda}\delta & \text{if $\delta\leq \sqrt{-2\lambda}$}, \\
\delta^2/2 & \text{if $\delta>\sqrt{-2\lambda}$}.
\end{cases}$$

\begin{lem}\label{lem-limsup}
Under the full condition of Theorem {\rm \ref{thm-main}}, for any $\delta>0$, 
$$\limsup_{t\rightarrow\infty}\frac{1}{t}\log Z_t^{\delta t}\leq -\Lambda_{\delta} 
\quad \text{${\mathbf P}_x$-a.s.}$$
\end{lem}
\pf \ 
For $t>0$, let $i$ be a nonnegative  integer such that $i\leq t<i+1$. 
Then for any $\delta>0$, 
\begin{equation}\label{eq-disc}
Z_t^{\delta t}\leq 
\sum_{k=1}^{Z_{i+1}}
{\bf 1}_{\left\{\sup_{i\leq s\leq i+1}|{\mathbf B}_s^{(i+1),k}|\geq \delta i\right\}}.
\end{equation}
The right hand side above is the number of particles alive at time $i+1$ such that 
their norms exceed $\delta i$ at least once during the time interval $[i,i+1]$. 
For $\varepsilon>0$, we define 
$$A_i:=
\left\{\sum_{k=1}^{Z_{i+1}}
{\bf 1}_{\left\{\sup_{i\leq s\leq i+1}|{\mathbf B}_s^{(i+1),k}|\geq \delta i\right\}}
>e^{(-\Lambda_{\delta}+\varepsilon)i}
\right\}  \quad (i\geq 0).$$
We will prove that there exist positive constants $c_1$, $c_2$ such that for all sufficiently large $i\geq 1$,
\begin{equation}\label{eq-upper}
{\mathbf P}_x(A_i)\leq c_1e^{-c_2i}.
\end{equation}
Since Borel-Cantelli's lemma implies that ${\mathbf P}_x(\text{$A_i$ i.o.})=0$, 
on an event of ${\mathbf P}_x$-full probability, 
there exists a natural number valued random variable $N\geq 1$ such that for all $i\geq N$,
$$
\sum_{k=1}^{Z_{i+1}}
{\bf 1}_{\left\{\sup_{i\leq s\leq i+1}|{\mathbf B}_s^{(i+1),k}|\geq \delta i\right\}}
\leq e^{(-\Lambda_{\delta}+\varepsilon)i}.
$$
Combining this with \eqref{eq-disc}, we see that  if $i\geq N$ and $i\leq t<i+1$, then 
$$
Z_t^{\delta t}
\leq e^{(-\Lambda_{\delta}+\varepsilon)i}
\leq \max\{1,e^{\Lambda_{\delta}-\varepsilon}\}e^{(-\Lambda_{\delta}+\varepsilon)t},
$$
which yields that 
$$\limsup_{t\rightarrow\infty}\frac{1}{t}\log Z_t^{\delta t}
\leq -\Lambda_{\delta}+\varepsilon.$$
By letting $\varepsilon\rightarrow+0$, the proof is complete.

Let us show \eqref{eq-upper}. 
We first see that by Chebyshev's inequality and the Markov property,
\begin{equation*}
\begin{split}
{\mathbf P}_x(A_i)
&\leq 
e^{(\Lambda_{\delta}-\varepsilon)i}
{\mathbf E}_x\left[
\sum_{k=1}^{Z_{i+1}}
{\bf 1}_{\left\{\sup_{i\leq s\leq i+1}|{\mathbf B}_s^{(i+1),k}|\geq \delta i\right\}}
\right]\\
&=e^{(\Lambda_{\delta}-\varepsilon)i}
{\mathbf E}_x\left[
\sum_{l=1}^{Z_i}
{\mathbf E}_{{\mathbf B}_i^{l}}
\left[\sum_{k=1}^{Z_1}
{\bf 1}_{\left\{\sup_{0\leq s\leq 1}|{\mathbf B}_s^{(1),k}|\geq \delta i\right\}}\right]\right].
\end{split}
\end{equation*}
Then by \eqref{eq-br-fk} and Lemma \ref{lem-exp}, 
the last expression above is equal to 
\begin{equation*}
\begin{split}
&e^{(\Lambda_{\delta}-\varepsilon)i}
E_x\left[e^{A_i^{(Q-1)\mu}}E_{B_i}\left[e^{A_1^{(Q-1)\mu}}; \sup_{0\leq s\leq 1}|B_s|\geq \delta i\right]\right]\\
&=
e^{(\Lambda_{\delta}-\varepsilon)i}
E_x\left[e^{A_{i+1}^{(Q-1)\mu}} ; \sup_{i\leq s\leq i+1}|B_s|\geq \delta i\right].
\end{split}
\end{equation*}

Under the condition that $\sup_{i\leq s\leq i+1}|B_s|\geq \delta i$, 
$|B_{i+1}|$ may be close to  $\delta i$ with high probability 
because the fluctuation of the Brownian motion is so small 
that \eqref{eq-fluct} below holds for any $\theta>0$. 
Taking this into account, for $\alpha\in (0,\delta)$, we let  
\begin{equation*}
\begin{split}
&E_x\left[e^{A_{i+1}^{(Q-1)\mu}} ; \sup_{i\leq s\leq i+1}|B_s|\geq \delta i\right]\\
&=E_x\left[e^{A_{i+1}^{(Q-1)\mu}} ; 
\sup_{i\leq s\leq i+1}|B_s|\geq \delta i, |B_{i+1}|\geq (\delta-\alpha)(i+1)\right]\\
&+E_x\left[e^{A_{i+1}^{(Q-1)\mu}} ; 
\sup_{i\leq s\leq i+1}|B_s|\geq \delta i, |B_{i+1}|< (\delta-\alpha)(i+1)\right]
={\rm (I)}+{\rm (II)}.
\end{split}
\end{equation*}
Then by Proposition \ref{prop-non-cpt}, 
$$
{\rm (I)}
\leq E_x\left[e^{A_{i+1}^{(Q-1)\mu}}; |B_{i+1}|\geq (\delta-\alpha)(i+1)\right]
\leq e^{(-\Lambda_{\delta-\alpha}+\varepsilon/2)(i+1)}
$$
for all sufficiently large $i\geq 1$. 
Therefore, if we take $\alpha>0$ so small that 
$$-C_1:=-\Lambda_{\delta-\alpha}+\Lambda_{\delta}-\frac{\varepsilon}{2}<0,$$
then 
$e^{(\Lambda_{\delta}-\varepsilon)i}{\rm (I)}\leq ce^{-C_1i}$ 
for all sufficiently large $i\geq 1$.

By the triangle inequality and the Markov property,
\begin{equation*}
\begin{split}
{\rm (II)}
&\leq E_x\left[e^{A_{i+1}^{(Q-1)\mu}}; 
\sup_{i\leq s\leq i+1}|B_s-B_{i+1}|>\delta i-|B_{i+1}|, 
|B_{i+1}|< (\delta-\alpha)(i+1)\right]\\
&\leq E_x
\left[
e^{A_{i+1}^{(Q-1)\mu}}; 
\sup_{i\leq s\leq i+1}|B_s-B_{i+1}|>\delta i-(\delta-\alpha)(i+1)\right]\\
&=E_x
\left[e^{A_i^{(Q-1)\mu}}
E_{B_i}\left[e^{A_1^{(Q-1)\mu}}; 
\sup_{0\leq s\leq 1}|B_s-B_1|\geq \alpha i-(\delta-\alpha) \right]
\right].
\end{split}
\end{equation*}
Then for any $x\in {\mathbb R}^d$ and $\theta>0$,
\begin{equation*}
\begin{split}
&E_x\left[e^{A_1^{(Q-1)\mu}}; \sup_{0\leq s\leq 1}|B_s-B_1|\geq \alpha i-(\delta-\alpha) \right]\\
&\leq e^{-\theta\{\alpha i-(\delta-\alpha)\}}
E_x\left[e^{A_1^{(Q-1)\mu}}e^{\theta \sup_{0\leq s\leq 1}|B_s-B_1|}\right]\\
&\leq e^{-\theta\{\alpha i-(\delta-\alpha)\}}E_x\left[e^{2A_1^{(Q-1)\mu}}\right]^{1/2}
E_x\left[e^{2\theta\sup_{0\leq s\leq 1}|B_s-B_1|}\right]^{1/2}.
\end{split}
\end{equation*}
Since 
$$\sup_{x\in {\mathbb R}^d}E_x\left[e^{2A_1^{(Q-1)\mu}}\right]<\infty$$
(see, e.g., \cite[Theorem 6.1]{ABM91}) and
\begin{equation}\label{eq-fluct}
E_x\left[e^{\theta\sup_{0\leq s\leq 1}|B_s-B_1|}\right]
=E_0\left[e^{\theta\sup_{0\leq s\leq 1}|B_s|}\right]<\infty,
\end{equation}
we have by \eqref{eq-comp-funct},
$${\rm (II)}
\leq 
c_1e^{-\theta\{\alpha i-(\delta-\alpha)\}}
E_x\left[e^{A_i^{(Q-1)\mu}}\right]
\leq c_2e^{-\theta\{\alpha i-(\delta-\alpha)\}}e^{-\lambda i}$$
for all sufficiently large $i\geq 1$. 
Therefore,
$$e^{(\Lambda_{\delta}-\varepsilon)i}{\rm (II)}
\leq c_2e^{(\Lambda_{\delta}-\varepsilon)i}
e^{-\theta\{\alpha i-(\delta-\alpha)\}}e^{-\lambda i}
=e^{\theta(\delta-\alpha)}
e^{(-\lambda -\theta\alpha+\Lambda_{\delta}-\varepsilon)i}.$$
Furthermore, by taking $\theta>0$ so large that 
$$-C_2:=-\lambda-\theta\alpha+\Lambda_{\delta}-\varepsilon<0,$$
we obtain 
$e^{(\Lambda_{\delta}-\varepsilon)i}{\rm (II)}
\leq c'e^{-C_2 i}$ for all sufficiently large $i\geq 1$, whence \eqref{eq-upper} holds. 
\qed
\medskip

\begin{lem}\label{lem-liminf}
Suppose that $\nu\in {\cal K}_{\infty}(1)$. 
If $\lambda<0$ and ${\mathbf P}_x(M_{\infty}>0)>0$, 
then for any $\delta\in (0,\sqrt{-\lambda/2})$, 
$$\liminf_{t\rightarrow\infty}\frac{1}{t}\log Z_t^{\delta t}
\geq -\lambda-\sqrt{-2\lambda}\delta \quad \text{${\mathbf P}_x(\cdot \mid M_{\infty}>0)$-a.s.}$$
\end{lem} 

To show this lemma, we follow the idea of  \cite[Lemma 4]{BH14}. 
Fix $p\in [0,1)$. 
Under the event that $M_{\infty}>0$, 
$Z_{tp}$ is approximately $e^{-\lambda pt}$ 
by Lemma \ref{lem-growth}. 
Moreover, each of the particles alive at time $tp$ will be on 
$\{x\in {\mathbb R}^d \mid |x|\geq \delta t\}$ at time $t$
with probability about $\exp\left(-\delta^2 t/\{2(1-p)\}\right)$ 
(see \eqref{eq-i} above).  
Therefore, $Z_t^{\delta t}$ has a rough lower bound  
\begin{equation}\label{eq-approx-lower}
\exp\left(-\lambda tp -\frac{\delta^2 t}{2(1-p)}\right).
\end{equation}
We will then choose $p$ attaining  the maximum of the right hand side 
as \eqref{eq-p_0} below.
\medskip

\pf \ 
Recall first that ${\mathbf B}_t^k$ is the position of the $k$th particle alive at time $t$. 
When this particle splits, we choose one of its children and follow its trajectory. 
By repeating this procedure inductively,  
we can construct a trajectory starting from ${\mathbf B}_t^k$. 
If we denote by ${\mathbf B}_s^{t,k}$ the position of such trajectory at time $s$ ($s\geq t$), 
then the law of $\{{\mathbf B}_s^{t,k}\}_{s\geq t}$ is identified 
with that of the Brownian motion on ${\mathbb R}^d$
starting at ${\mathbf B}_t^k$.

Fix $p\in [0,1)$. Then by the triangle inequality,
\begin{equation}\label{eq-event}
\begin{split}
&\left\{\text{$|{\mathbf B}_s^{np, k}|>\delta s$ for all $s\in [n,n+1]$}\right\}\\
&\supset \left\{
|{\mathbf B}_n^{np, k}|>|{\mathbf B}_n^{np, k}-{\mathbf B}_{np}^{np,k}|>\delta(n+1)+1, 
\sup_{n\leq s\leq n+1}|{\mathbf B}_s^{np, k}-{\mathbf B}_n^{np, k}|<1\right\}
=:E_n^{k}.
\end{split}
\end{equation}
For any $\alpha>0$, we see by the Markov property that 
\begin{equation}\label{eq-e}
\begin{split}
&{\mathbf P}_x\left(
\sum_{k=1}^{Z_{np}}
{\bf 1}_{E_n^k}\leq e^{(-\Lambda_{\delta}-\varepsilon)n}, 
Z_{np}\geq e^{(-\lambda p-\alpha)n}\right)\\
&\leq 
{\mathbf P}_x\left(
\sum_{k=1}^{[e^{(-\lambda p-\alpha)n}]}{\bf 1}_{E_n^k}
\leq e^{(-\Lambda_{\delta}-\varepsilon)n}, 
Z_{np}\geq e^{(-\lambda p-\alpha)n}\right)\\
&={\mathbf E}_x\left[{\mathbf P}_{{\mathbf B}_{np}}
\left(\sum_{k=1}^{[e^{(-\lambda p-\alpha)n}]}
{\bf 1}_{F_n^k}\leq e^{(-\Lambda_{\delta}-\varepsilon)n}\right)
; Z_{np}\geq e^{(-\lambda p-\alpha)n}\right]
\end{split}
\end{equation}
for 
$$F_n^k:=\left\{
|{\mathbf B}_{n-np}^{0,k}|\geq |{\mathbf B}_{n-np}^{0,k}-{\mathbf B}_0^{0,k}|>\delta(n+1)+1, 
\sup_{n-np\leq s\leq n-np+1}|{\mathbf B}_s^{0,k}-{\mathbf B}_{n-np}^{0,k}|<1\right\}.$$
Here we let ${\bf x}=(x^1,\dots, x^m)\in {\mathbf X}$ for $m\geq [e^{(-\lambda p-\alpha)n}]$. 
Since the events $F_n^k$ ($1\leq k\leq m$) are independent under ${\mathbf P}_{{\bf x}}$, 
it follows by Chebyshev's inequality that 
\begin{equation}\label{eq-enk}
\begin{split}
&{\mathbf P}_{{\mathbf x}}
\left(\sum_{k=1}^{[e^{(-\lambda p-\alpha)n}]}
{\bf 1}_{F_n^k}\leq e^{(-\Lambda_{\delta}-\varepsilon)n}\right)
={\mathbf P}_{{\mathbf x}}\left(\exp\left(-\sum_{k=1}^{[e^{(-\lambda p-\alpha)n}]}
{\bf 1}_{F_n^k}\right)
\geq e^{-e^{(-\Lambda_{\delta}-\varepsilon)n}}\right)\\
&\leq 
e^{e^{(-\Lambda_{\delta}-\varepsilon)n}}
{\mathbf E}_{\mathbf x}
\left[\exp\left(-\sum_{k=1}^{[e^{(-\lambda p-\alpha)n}]}{\bf 1}_{F_n^k}\right)\right]
=e^{e^{(-\Lambda_{\delta}-\varepsilon)n}}
\prod_{k=1}^{[e^{(-\lambda p-\alpha)n}]}
{\mathbf E}_{x^k}\left[\exp\left(-{\bf 1}_{F_n^k}\right)\right]\\
&=e^{e^{(-\Lambda_{\delta}-\varepsilon)n}}
\prod_{k=1}^{[e^{(-\lambda p-\alpha)n}]}
\left\{1-(1-e^{-1}){\mathbf P}_{x^k}(F_n^k)\right\}.
\end{split}
\end{equation}

Note that 
\begin{equation*}
\begin{split}
{\mathbf P}_{x}(F_n^k)
&=P_x\left(|B_{n-np}|\geq |B_{n-np}-x|>\delta(n+1)+1, 
\sup_{n-np\leq s\leq n-np+1}|B_s-B_{n-np}|<1\right)\\
&=P_x\left(|B_{n-np}|\geq |B_{n-np}-x|>\delta(n+1)+1\right) 
P_0\left(\sup_{0\leq s\leq 1}|B_s|<1\right).
\end{split}
\end{equation*}
Then
\begin{equation}\label{eq-cov}
\begin{split}
&P_x\left(|B_{n-np}|\geq |B_{n-np}-x|>\delta(n+1)+1\right) \\
&=\frac{1}{(2\pi(n-np))^{d/2}}\int_{|y|\geq |y-x|>\delta(n+1)+1}
\exp\left(-\frac{|y-x|^2}{2(n-np)}\right)\,{\rm d}y\\
&=\frac{1}{(2\pi(n-np))^{d/2}}\int_{|w-x|\geq |w|>\delta(n+1)+1}
\exp\left(-\frac{|w|^2}{2(n-np)}\right)\,{\rm d}w.
\end{split}
\end{equation}
Let $\langle \cdot, \cdot \rangle$ be the standard inner product on ${\mathbb R}^d$. 
Then the last term of \eqref{eq-cov} is not less than 
\begin{equation*}
\begin{split}
&\frac{1}{(2\pi(n-np))^{d/2}}\int_{|w|>\delta(n+1)+1, \langle w,x\rangle<0}
\exp\left(-\frac{|w|^2}{2(n-np)}\right)\,{\rm d}w\\
&=\frac{c}{(2\pi(n-np))^{d/2}}
\int_{\delta(n+1)+1}^{\infty}\exp\left(-\frac{r^2}{2(n-np)}\right)r^{d-1}\,{\rm d}r.
\end{split}
\end{equation*}
Furthermore, by the change of variables and \eqref{eq-i}, 
the last expression above is equal to 
\begin{equation*}
\begin{split}
\frac{c}{(2\pi)^{d/2}}\int_{(\delta(n+1)+1)/\sqrt{n-np}}^{\infty}e^{-r^2/2}r^{d-1}\,{\rm d}r
&\asymp \left(\frac{\delta(n+1)+1}{\sqrt{n-np}}\right)^{d-2}
\exp\left(-\frac{(\delta(n+1)+1)^2}{2(n-np)}\right)\\
&\asymp c_{\delta,p} n^{(d-2)/2}\exp\left(-\frac{\delta^2n}{2(1-p)}\right) \quad (n\rightarrow\infty).
\end{split}
\end{equation*}
Since this shows that 
$$
{\mathbf P}_x(F_n^k)
\geq c'_{\delta,p}n^{(d-2)/2}\exp\left(-\frac{\delta^2n}{2(1-p)}\right)
=:q_{n,p},
$$
the inequality $1-x\leq e^{-x}$ implies that the last term in \eqref{eq-enk} is less than 
\begin{equation}\label{eq-upper-enk}
e^{e^{(-\Lambda_{\delta}-\varepsilon)n}}
(1-(1-e^{-1})q_{n,p})^{[e^{(-\lambda p-\alpha)n}]}
\leq 
e^{e^{(-\Lambda_{\delta}-\varepsilon)n}}
\exp\left(-(1-e^{-1})q_{n,p}[e^{(-\lambda p-\alpha)n}]\right).
\end{equation}
Here we observe that 
\begin{equation}\label{eq-q-asymp}
\begin{split}
q_{n,p}[e^{(-\lambda p-\alpha)n}]
&\asymp 
n^{(d-2)/2}\exp\left(-\frac{\delta^2n}{2(1-p)}\right)
e^{(-\lambda p-\alpha)n}\\
&=n^{(d-2)/2}e^{-\alpha n}
\exp\left\{\left(-\lambda p-\frac{\delta^2}{2(1-p)}\right)n\right\} 
\quad (n\rightarrow\infty).
\end{split}
\end{equation}

Let 
$$F(p):=-\lambda p-\frac{\delta^2}{2(1-p)} \quad (0\leq p<1).$$
Then $F(p)$ has the maximal value
$F(p_0)=-\Lambda_{\delta}$ for
\begin{equation}\label{eq-p_0}
p_0=
\begin{cases}
0 & \text{if $\delta\geq \sqrt{-2\lambda}$,}\\
1-\delta/{\sqrt{-2\lambda}}\in (0,1) & \text{if $\delta<\sqrt{-2\lambda}$.}
\end{cases}
\end{equation}
In what follows, we assume that  $p=p_0$. Then by \eqref{eq-q-asymp},
\begin{equation*}
\begin{split}
e^{e^{(-\Lambda_{\delta}-\varepsilon)n}}
\exp\left(-(1-e^{-1})q_{n,p}[e^{(-\lambda p-\alpha)n}]\right)
&\leq \exp\left\{e^{(-\Lambda_{\delta}-\varepsilon)n}-cn^{(d-2)/2}
e^{(-\Lambda_{\delta}-\alpha) n}\right\}\\
&=\exp\left\{-e^{(-\Lambda_{\delta}-\alpha) n}\left(cn^{(d-2)/2}-e^{-(\varepsilon-\alpha)n}\right)\right\}.
\end{split}
\end{equation*}
In particular, if we further assume that $\delta\in (0,\sqrt{-\lambda/2})$ and $\alpha<\varepsilon$, 
then $\Lambda_{\delta}<0$ and thus the right hand side above converges exponentially to $0$ as $n\rightarrow\infty$. 
Therefore, by \eqref{eq-e}, \eqref{eq-enk} and \eqref{eq-upper-enk}, 
$$\sum_{n=1}^{\infty}{\mathbf P}_x\left(
\sum_{k=1}^{Z_{np}}
{\bf 1}_{E_n^k}
\leq e^{(-\Lambda_{\delta}-\varepsilon)n}, 
Z_{np}\geq e^{(-\lambda p-\alpha)n}\right)<\infty.$$
By Borel-Cantelli's lemma,  if we set 
$$G_n:=\left\{\sum_{k=1}^{Z_{np}}{\bf 1}_{E_n^k}
>e^{(-\Lambda_{\delta}-\varepsilon)n}\right\}
\bigcup
\left\{Z_{np}< e^{(-\lambda p-\alpha)n}\right\},$$
then ${\mathbf P}_x(\text{$G_n$ eventually})=1$. 
Since
$${\mathbf P}_x\left(\text{$Z_{np}\geq e^{(-\lambda p-\alpha)n}$ eventually}\mid M_{\infty}>0\right)=1$$
by Lemma \ref{lem-growth}, we obtain
$$
{\mathbf P}_x\left(
\text{$\sum_{k=1}^{Z_{np}}{\bf 1}_{E_n^k}
>e^{(-\Lambda_{\delta}-\varepsilon)n}$ eventually} \mid M_{\infty}>0\right)=1.
$$
Hence under the event $\{M_{\infty}>0\}$, 
it follows by \eqref{eq-event} that for all sufficiently large $t>0$,
$$
Z_t^{\delta t}
=\sum_{k=1}^{Z_t}
{\bf 1}_{\{|{\mathbf B}_t^k|>\delta t\}}
\geq \sum_{k=1}^{Z_{[t]p}}
{\bf 1}_{\{\text{$|{\mathbf B}_s^{[t]p, k}|>\delta s$ 
for all $s\in [[t],[t]+1]$}\}}
\geq \sum_{k=1}^{Z_{[t]p}}
{\bf 1}_{E_{[t]}^{k}}\geq e^{(-\Lambda_{\delta}-\varepsilon)[t]},
$$
which yields that
$$\liminf_{t\rightarrow\infty}\frac{1}{t}\log Z_t^{\delta t}\geq -\Lambda_{\delta}-\varepsilon.$$
We thus complete the proof 
by letting first $\alpha\rightarrow +0$ and then $\varepsilon\rightarrow +0$. 
\qed
\medskip

{\it Proof of Theorem {\rm \ref{thm-main}}}. 
The assertion (i) follows from Lemma \ref{lem-limsup} 
because $\Lambda_{\delta}>0$ for any $\delta>\sqrt{-\lambda/2}$ and 
$Z_t^{\delta t}$ takes a nonnegative integer value for any $t>0$.
The assertion (ii) is a consequence of Lemmas \ref{lem-growth}, \ref{lem-limsup} and \ref{lem-liminf}. 
\qed
\medskip

{\it Proof of Theorem {\rm \ref{thm-direction}}}.
Let $r$  be a unit vector in ${\mathbb R}^d$. 
We first prove (i). 
Let $B^{r,+}:=\{x\in {\mathbb R}^d \mid \langle x, r\rangle\geq 0\}$ and $Z_t^{r,+}:=Z_t(B^{r,+})$. 
Then by \eqref{eq-limit-thm} and Lemma \ref{lem-growth},  
\begin{equation}\label{eq-growth-p}
\lim_{t\rightarrow\infty}\frac{1}{t}\log Z_t^{r,+}=-\lambda, 
\quad \text{${\mathbf P}_x(\cdot \mid M_{\infty}>0)$-a.s.}
\end{equation}
Therefore, (i) holds for $\delta=0$. 

In what follows, we assume that $\delta>0$. 
Since $Z_t^{\delta t,r}\leq Z_t^{\delta t}$, we have by Lemma \ref{lem-limsup},
$$\limsup_{t\rightarrow\infty}\frac{1}{t}\log Z_t^{\delta t,r}
\leq \limsup_{t\rightarrow\infty}\frac{1}{t}\log Z_t^{\delta t}
\leq-\Lambda_{\delta}, \quad \text{${\mathbf P}_x$-a.s.}$$
For a fixed $p\in [0,1)$, if $\langle {\mathbf B}_{np}^{np,k},r\rangle\geq 0$, then 
\begin{equation*}
\begin{split}
&\left\{\text{$\langle {\mathbf B}_s^{np, k},r\rangle>\delta s$ for all $s\in [n,n+1]$}\right\}\\
&\supset \left\{
\langle {\mathbf B}_n^{np, k}-{\mathbf B}_{np}^{np,k},r\rangle>\delta(n+1)+1, 
\inf_{n\leq s\leq n+1}\langle {\mathbf B}_s^{np, k}-{\mathbf B}_n^{np, k},r\rangle>-1\right\}
=:E_n^{r,k}.
\end{split}
\end{equation*}
Hence by following the argument of \cite[Lemma 4]{BH14} or Lemma \ref{lem-liminf}, 
we see from \eqref{eq-growth-p} that if $\delta\in (0,\sqrt{-\lambda/2})$, then 
$$
\sum_{n=1}^{\infty}{\mathbf P}_x\left(
\sum_{1\leq k\leq Z_{np}, \langle {\mathbf B}_{np}^{np,k},r\rangle\geq 0}
{\bf 1}_{E_n^{r,k}}
\leq e^{(-\Lambda_{\delta}-\varepsilon)n}, 
Z_{np}^{r,+}\geq e^{(-\lambda p-\alpha)n}\right)<\infty$$
and 
$$\liminf_{t\rightarrow\infty}\frac{1}{t}\log Z_t^{\delta t,r}
\geq -\lambda-\sqrt{-2\lambda}\delta, \quad \text{${\mathbf P}_x(\cdot \mid M_{\infty}>0)$-a.s.}$$
We thus arrive at  (i). 

We next prove (ii). 
By the same argument as in Corollary \ref{cor-main}, we get \eqref{eq-direction-0}. 
To show \eqref{eq-direction}, we note that if $\{r_1,\dots, r_d\}$ is an orthonormal basis of ${\mathbb R}^d$ such that $r_1=r$, 
then ${\mathbf B}_t^{K_r(t)}=\sum_{j=1}^d\langle {\mathbf B}_t^{K_r(t)},r_j\rangle r_j$ and 
$$|{\mathbf B}_t^{K_r(t)}|^2=\sum_{j=1}^d\langle {\mathbf B}_t^{K_r(t)},r_j\rangle^2
=\langle {\mathbf B}_t^{K_r(t)},r\rangle^2+\sum_{j=2}^d\langle {\mathbf B}_t^{K_r(t)},r_j\rangle^2.$$
Hence  by \eqref{eq-direction-0} and Corollary \ref{cor-main}, we have ${\mathbf P}_x(\cdot \mid M_{\infty}>0)$-a.s.,
$$\frac{1}{t^2}\sum_{j=2}^d\langle {\mathbf B}_t^{K_r(t)},r_j\rangle^2
=\frac{1}{t^2}|{\mathbf B}_t^{K_r(t)}|^2-\frac{1}{t^2}\langle {\mathbf B}_t^{K_r(t)},r\rangle^2
\leq \left(\frac{L_t}{t}\right)^2-\left(\frac{L_t^r}{t}\right)^2\rightarrow 0 \quad (t\rightarrow\infty),$$
which yields \eqref{eq-direction}.
\qed
\medskip

Related to Corollary \ref{cor-main}, we finally reveal the decay rate of  
${\mathbf P}_x(L_t\geq \delta t)$ as $t\rightarrow\infty$ for $\delta>\sqrt{-\lambda/2}$. 
\begin{lem}\label{lem-positive}
Under the same setting as in Theorem {\rm \ref{thm-main}}, 
if ${\mathbf P}_x(M_{\infty}>0)>0$, then for any $\delta>\sqrt{-\lambda/2}$,
$$\lim_{t\rightarrow \infty}\frac{1}{t}\log {\mathbf P}_x\left(\frac{L_t}{t}\geq \delta\right)=\lim_{t\rightarrow \infty}\frac{1}{t}\log {\mathbf P}_x(Z_t^{\delta t}\geq 1)=-\Lambda_{\delta}.$$ 
\end{lem}

To give a lower bound of ${\mathbf P}_x(Z_t^{\delta t}\geq 1)$, 
we follow the same idea as Lemma \ref{lem-liminf}. 
For $p\in [0,1)$, 
$$
{\mathbf P}_x(Z_t^{\delta t}\geq 1)=1-{\mathbf P}_x(Z_t^{\delta t}=0)
={\mathbf E}_x\left[1-{\mathbf P}_{{\mathbf B}_{pt}}(Z_{(1-p)t}^{\delta t}=0)\right].
$$
By ignoring the starting point and the effect of branching, 
$1-{\mathbf P}_{{\mathbf B}_{pt}}(Z_{(1-p)t}^{\delta t}=0)$ has a rough lower bound
$$1-P_0(|B_{(1-p)t}|<\delta t)^{Z_{pt}}=1-\left(1-P_0(|B_{(1-p)t}|\geq \delta t)\right)^{Z_{pt}}
\sim Z_{pt}P_0(|B_{(1-p)t}|\geq \delta t) \quad (t\rightarrow\infty).$$
Furthermore, under the event that $M_{\infty}>0$, 
the right hand side above is approximately \eqref{eq-approx-lower} by Lemma \ref{lem-growth}.  
We can thus give  a rough lower bound of ${\mathbf P}_x(Z_t^{\delta t}\geq 1)$ by 
$$\exp\left(-\lambda pt-\frac{\delta^2 t}{2(1-p)}\right){\mathbf P}_x(M_{\infty}>0).$$
We will then take $p=p_0$ as in \eqref{eq-p_0} to maximize the right hand side above.
\medskip

\pf \ 
The first equality follows from the fact that 
${\mathbf P}_x(L_t\geq \delta t)={\mathbf P}_x(Z_t^{\delta t}\geq 1)$.

Let us discuss the second equality. Since 
$${\mathbf P}_x(Z_t^{\delta t}\geq 1)\leq {\mathbf E}_x[Z_t^{\delta t}]
=E_x\left[e^{A_t^{(Q-1)\mu}}; |B_t|\geq \delta t\right]$$
by the Chebyshev inequality, we have by Proposition \ref{prop-non-cpt},
$$\limsup_{t\rightarrow \infty}\frac{1}{t}\log {\mathbf P}_x(Z_t^{\delta t}\geq 1)\leq -\Lambda_{\delta}.$$ 

In what follows, we use the same notations as those in the proof of Lemma \ref{lem-liminf}. 
For any $t>0$, $p\in [0,1)$ and $\alpha>0$,
\begin{equation}\label{eq-l-1}
\begin{split}
{\mathbf P}_x(Z_t^{\delta t}\geq 1)
&={\mathbf P}_x\left(\bigcup_{k=1}^{Z_{tp}}\left\{|{\mathbf B}_t^{tp,k}|\geq \delta t\right\}\right)\\
&\geq {\mathbf P}_x\left(\bigcup_{k=1}^{Z_{tp}}\left\{|{\mathbf B}_t^{tp,k}|\geq \delta t\right\}, 
Z_{tp}\geq e^{(-\lambda p-\alpha)t}\right)\\
&={\mathbf E}_x\left[{\mathbf P}_{{\mathbf B}_{tp}}
\left(\bigcup_{k=1}^{m}\left\{|{\mathbf B}_{t-tp}^{0,k}|\geq \delta t\right\}\right)\Big|_{m=Z_{tp}};
Z_{tp}\geq e^{(-\lambda p-\alpha)t}\right].
\end{split}
\end{equation}
At the last equality, we used the Markov property. 

By \eqref{eq-lower-ball} and \eqref{eq-int-comp}, 
there exist $T>0$ and $c>0$ such that for any ${\bf x}=(x^1,\dots, x^m)\in {\bf X}$,  and $t\geq T$,
\begin{equation*}
\begin{split}
{\mathbf P}_{\bf x}\left(|{\mathbf B}_{t-tp}^{0,k}|\geq \delta t\right)
=P_{x^k}\left(|B_{t-tp}|\geq \delta t\right)
&\geq P_0\left(|B_{t-tp}|\geq \delta t\right)\\
&\geq c\exp\left(-\frac{\delta^2 t}{2(1-p)}\right)t^{(d-2)/2}=:f(t).
\end{split}
\end{equation*}
Since the events $\{|{\mathbf B}_{t-tp}^{0,k}|\geq \delta t\} \ (1\leq k\leq m)$ are independent under ${\mathbf P}_{\bf x}$,
we have for any $t\geq T$,
\begin{equation*}
\begin{split}
{\mathbf P}_{{\bf x}}
\left(\bigcup_{k=1}^{m}\left\{|{\mathbf B}_{t-tp}^{0,k}|\geq \delta t\right\}\right)
&=1-\prod_{k=1}^{m}
\left(1-{\mathbf P}_{x^k}
\left(|{\mathbf B}_{t-tp}^{0,k}|\geq \delta t\right)\right)\geq 1-(1-f(t))^m.
\end{split}
\end{equation*}
Hence if $t\geq T$, then the last term of \eqref{eq-l-1} is not less than 
\begin{equation}\label{eq-l-2}
{\mathbf E}_x\left[1-(1-f(t))^{Z_{tp}}; 
Z_{tp}\geq e^{(-\lambda p-\alpha)t}\right]
\geq \left\{1-(1-f(t))^{g(t)}\right\}
{\mathbf P}_x\left(Z_{tp}\geq e^{(-\lambda p-\alpha)t}\right)
\end{equation}
for $g(t):=e^{(-\lambda p-\alpha)t}$. 

If we take $\delta>\sqrt{-\lambda/2}$ and $p=p_0$ as in \eqref{eq-p_0}, then $\Lambda_{\delta}>0$ and thus 
$$f(t)g(t)=c\exp\left(-\frac{\delta^2 t}{2(1-p)}-(\lambda p+\alpha)t\right)t^{(d-2)/2}
=ce^{-(\Lambda_{\delta}+\alpha)t}t^{(d-2)/2}\rightarrow 0 \quad (t\rightarrow \infty).$$
By this together with elementary calculation, we obtain  
$$\liminf_{t\rightarrow\infty}\frac{1}{t}\log (1-(1-f(t))^{g(t)})\geq -\Lambda_{\delta}-\alpha.$$
We also see by Lemma \ref{lem-growth} that 
$${\mathbf P}_x\left(Z_{tp}\geq e^{(-\lambda p-\alpha)t}\right)
\geq {\mathbf P}_x\left(Z_{tp}\geq e^{(-\lambda p-\alpha)t}, M_{\infty}>0\right)
\rightarrow {\mathbf P}_x(M_{\infty}>0) \quad (t\rightarrow\infty).$$
Combining these two inequalities above with \eqref{eq-l-1} and \eqref{eq-l-2}, we get
$$\liminf_{t\rightarrow\infty}\frac{1}{t}\log {\mathbf P}_x(Z_t^{\delta t}\geq 1)\geq -\Lambda_{\delta}-\alpha.$$
Letting $\alpha\rightarrow +0$, we arrive at our assertion.
\qed

\begin{rem}\rm
Let $r$ be a unit vector in ${\mathbb R}^d$. 
Then Lemma \ref{lem-positive} is valid by replacing 
$L_t$ and $Z_t^{\delta t}$ with $L_t^r$ and $Z_t^{\delta t,r}$, respectively. 
We can show this by the same way as for the proof of Lemma \ref{lem-positive}.
\end{rem}

\appendix
\section{Appendix}\label{appendix-1}
\subsection{Proof of \eqref{eq-lower-ball} and \eqref{eq-lower-ball-0}}\label{appendix-a-1}
In this subsection, we prove that for any $R>0$, 
$$
P_x(|B_t|\geq R)\geq P_0(|B_t|\geq R) \quad \text{for any $x\in {\mathbb R}^d$ and  $t>0$}.
$$
This inequality implies \eqref{eq-lower-ball} and \eqref{eq-lower-ball-0}.

For $x\ne 0$, 
\begin{equation}\label{eq-sph}
\begin{split}
\int_{|y|\geq R}\exp\left(-\frac{|y-x|^2}{2t}\right)\,{\rm d}y
&=\exp\left(-\frac{|x|^2}{2t}\right) 
\int_{|y|\geq R}\exp\left(-\frac{|y|^2}{2t}+\frac{\langle x/|x|,y\rangle}{t}|x|\right)\,{\rm d}y\\
&=\exp\left(-\frac{|x|^2}{2t}\right) 
\int_{|y|\geq R}\exp\left(-\frac{|y|^2}{2t}+\frac{y_1}{t}|x|\right)\,{\rm d}y
\end{split}
\end{equation}
by the spherical symmetry. 
Define for $s\geq 0$,
$$F(s):=\exp\left(-\frac{s^2}{2t}\right) 
\int_{|y|\geq R}\exp\left(-\frac{|y|^2}{2t}+\frac{s y_1}{t}\right)\,{\rm d}y.$$
If we set 
$$|y|_d=\sqrt{y_2^2+\dots+y_d^2}, \quad K=K(y_2,\dots, y_d)=\sqrt{R^2-|y|_d^2},$$
then 
\begin{equation*}
\begin{split}
&F'(s)
=\exp\left(-\frac{s^2}{2t}\right)\int_{|y|\geq R}
\exp\left(-\frac{|y|_d^2}{2t}\right)\frac{y_1-s}{t}\exp\left(\frac{2s y_1-y_1^2}{2t}\right)\,{\rm d}y\\
&=\exp\left(-\frac{s^2}{2t}\right)\int_{|y|_d\geq R}\exp\left(-\frac{|y|_d^2}{2t}\right)
\left(\int_{{\mathbb R}}\frac{y_1-s}{t}\exp\left(\frac{2sy_1-y_1^2}{2t}\right)\,{\rm d}y_1\right)\,{\rm d}y_2\dots{\rm d}y_d\\
&+\exp\left(-\frac{s^2}{2t}\right)\int_{|y|_d<R}\exp\left(-\frac{|y|_d^2}{2t}\right)
\left(\int_{|y_1|\geq K}\frac{y_1-s}{t}\exp\left(\frac{2sy_1-y_1^2}{2t}\right)\,{\rm d}y_1\right)\,{\rm d}y_2\dots{\rm d}y_d\\
&=\exp\left(-\frac{s^2+R^2}{2t}\right)\int_{|y|_d<R}
\left(\exp\left(\frac{sK}{t}\right)-\exp\left(-\frac{sK}{t}\right)\right)
\,{\rm d}y_2\dots{\rm d}y_d\geq 0.
\end{split}
\end{equation*}
Therefore, 
$$P_x(|B_t|\geq R)=\frac{1}{(2\pi t)^{d/2}}F(|x|)\geq \frac{1}{(2\pi t)^{d/2}}F(0)=P_x(|B_t|\geq R),$$
which is our assertion.

\subsection{Proof of \eqref{eq-iv}}\label{appendix-a-2}

In this subsection, we show \eqref{eq-iv}.  
Here we use the same notations as in the proof of Proposition \ref{prop-asymp}.
We first assume that $\delta<\sqrt{-2\lambda}$.  
Then for all sufficiently large $t$, 
$$S=\left(\sqrt{-\lambda}-\frac{\delta}{\sqrt{2}}-\frac{\delta}{\sqrt{2}}\frac{T}{t-T}\right)\sqrt{t-T}>0.$$
We let 
\begin{equation}\label{eq-iv-split}
\begin{split}
({\rm IV})&=\int_{-\infty}^0e^{-v^2}
\left(\frac{2\sqrt{-\lambda}}{v+\sqrt{v^2+2\sqrt{-2\lambda}\delta t}}\right)^d\frac{1}{\sqrt{v^2+2\sqrt{-2\lambda}\delta t}}\,{\rm d}v\\
&+\int_0^S e^{-v^2}
\left(\frac{2\sqrt{-\lambda}}{v+\sqrt{v^2+2\sqrt{-2\lambda}\delta t}}\right)^d\frac{1}{\sqrt{v^2+2\sqrt{-2\lambda}\delta t}}\,{\rm d}v\\
&=({\rm IV})_1+({\rm IV})_2.
\end{split}
\end{equation}
If $v\leq 0$, then 
\begin{equation}\label{eq-calcu}
\frac{1}{v+\sqrt{v^2+2\sqrt{-2\lambda}\delta t}}
=\frac{\sqrt{v^2+2\sqrt{-2\lambda}\delta t}-v}{2\sqrt{-2\lambda}\delta t}
\asymp \frac{\sqrt{v^2+2\sqrt{-2\lambda}\delta t}}{t} \quad (t\rightarrow\infty)
\end{equation}
so that we get 
\begin{equation*}
\begin{split}
({\rm IV})_1
&\asymp 
\int_{-\infty}^0e^{-v^2}
\left(\frac{\sqrt{v^2+2\sqrt{-2\lambda}\delta t}}{t}\right)^d
\frac{1}{\sqrt{v^2+2\sqrt{-2\lambda}\delta t}}\,{\rm d}v\\
&=\frac{1}{t^d}\int_0^{\infty}e^{-v^2}\left(v^2+2\sqrt{-2\lambda}\delta t\right)^{(d-1)/2}\,{\rm d}v 
\quad (t\rightarrow\infty).
\end{split}
\end{equation*}
Since there exist positive constants $c_1$, $c_2$ such that 
\begin{equation}\label{eq-poly}
c_1t^{(d-1)/2}
\leq \left(v^2+2\sqrt{-2\lambda}\delta t\right)^{(d-1)/2}
\leq c_2(v^{d-1}+t^{(d-1)/2}) \quad \text{for all $v\geq 0$},
\end{equation}
we have 
\begin{equation}\label{eq-int-iv}
\begin{split}
\frac{c_1}{t^{(d+1)/2}}\int_0^{\infty}e^{-v^2}\,{\rm d}v
&\leq \frac{1}{t^d}\int_0^{\infty}e^{-v^2}\left(v^2+2\sqrt{-2\lambda}\delta t\right)^{(d-1)/2}\,{\rm d}v\\
&\leq c_2\left(\frac{1}{t^d}\int_0^{\infty}
e^{-v^2}v^{d-1}\,{\rm d}v
+\frac{1}{t^{(d+1)/2}}\int_0^{\infty}e^{-v^2}\,{\rm d}v\right),
\end{split}
\end{equation}
which implies that 
\begin{equation}\label{eq-a-iv}
({\rm IV})_1\asymp t^{-(d+1)/2} \quad (t\rightarrow\infty).
\end{equation}
On the other hand, there exist positive constants $c_3$, $c_4$ 
such that 
$$
c_3 t^{-1/2}
\leq \frac{1}{v+\sqrt{v^2+2\sqrt{-2\lambda}\delta t}}
\leq \frac{1}{\sqrt{v^2+2\sqrt{-2\lambda}\delta t}}
\leq c_4 t^{-1/2} \quad (0\leq v\leq S) 
$$
because
$$S\leq \left(\sqrt{-\lambda}-\frac{\delta}{\sqrt{2}}\right)\sqrt{t}.$$
Noting that $S\rightarrow\infty$ as $t\rightarrow\infty$, we have 
$$
({\rm IV})_2
\asymp t^{-(d+1)/2} \quad (t\rightarrow\infty).
$$
Combining this with \eqref{eq-a-iv} and \eqref{eq-iv-split}, we get   
$$({\rm IV})=({\rm IV})_1+({\rm IV})_2\asymp t^{-(d+1)/2} \quad (t\rightarrow\infty).$$

We next assume that $\delta=\sqrt{-2\lambda}$. 
Since
$$S=-\frac{\delta}{\sqrt{2}}\frac{T}{\sqrt{t-T}}\rightarrow0 \quad (t\rightarrow\infty),$$
we have by \eqref{eq-poly} and \eqref{eq-int-iv}, 
\begin{equation*}
\begin{split}
({\rm IV})
\asymp 
\frac{1}{t^d}\int_{-\infty}^0
e^{-v^2}
(v^2+2\sqrt{-2\lambda}\delta t)^{(d-1)/2}\,{\rm d}v
&=\frac{1}{t^d}
\int_0^{\infty}
e^{-v^2}
(v^2+2\sqrt{-2\lambda}\delta t)^{(d-1)/2}\,{\rm d}v\\
&\asymp \frac{1}{t^{(d+1)/2}} \quad (t\rightarrow\infty).
\end{split}
\end{equation*}

We finally assume that $\delta>\sqrt{-2\lambda}$. 
Since 
\begin{equation}\label{eq-s-comp} 
S
=\left(\sqrt{-\lambda}-\frac{\delta}{\sqrt{2}}-\frac{\delta}{\sqrt{2}}\frac{T}{t-T}\right)\sqrt{t-T}
<0, 
\end{equation}
we obtain by \eqref{eq-calcu},
\begin{equation}\label{eq-int}
\begin{split}
({\rm IV})
&\asymp 
\frac{1}{t^d}\int_{-\infty}^S
e^{-v^2}
(v^2+2\sqrt{-2\lambda}\delta t)^{(d-1)/2}\,{\rm d}v\\
&=\frac{1}{t^d}
\int_{-S}^{\infty}
e^{-v^2}
(v^2+2\sqrt{-2\lambda}\delta t)^{(d-1)/2}\,{\rm d}v \quad (t\rightarrow\infty).
\end{split}
\end{equation}
Note that by \eqref{eq-s-comp}, 
$$v^2\leq v^2+2\sqrt{-2\lambda}\delta t\leq c v^2 \quad \text{for all $v\geq -S$}$$
for some $c>0$. Since $S\asymp -\sqrt{t}$ as $t\rightarrow\infty$, we have by \eqref{eq-int-comp},
\begin{equation*}
\begin{split}
({\rm IV})
&\asymp 
\frac{1}{t^d}\int_{-S}^{\infty}
e^{-v^2}v^{d-1}\,{\rm d}v
\sim \frac{1}{t^d}e^{-S^2} \left(-S\right)^{d-2}
\asymp  e^{-(-\lambda-\sqrt{-2\lambda}\delta) t}e^{-\delta^2t/2}t^{-(d+2)/2}
\quad (t\rightarrow\infty).
\end{split}
\end{equation*}
Therefore, we have proved \eqref{eq-iv}.
\bigskip

\noindent
{\bf Acknowledgements} \ 
The author would like to thank Professor Masayoshi Takeda 
for his valuable comments on the draft of this paper. 
He is grateful to the two referees for their careful reading of the manuscript  
and constructive suggestions. 
In particular, they pointed out a mistake in Lemma \ref{lem-liminf} 
and one of them suggested showing Lemma \ref{lem-positive}. 
Thanks are also due to the editors for their helpful suggestions 
on the improvement of Introduction.

\address{
Yuichi Shiozawa\\
Department of Mathematics\\
Graduate School of Science\\ 
Osaka University\\
Toyonaka, Osaka, 560-0043,
Japan
}
{\texttt{shiozawa@math.sci.osaka-u.ac.jp}}
\end{document}